\documentclass[12pt]{amsart}

\usepackage{amsmath,amsfonts,amsthm}
\usepackage{tikz}

\renewcommand{\theequation}{\thesection.\arabic{equation}}
\newtheorem{theorem}{Theorem}
\newtheorem{lemma}{Lemma}
\newtheorem{proposition}{Proposition}
\newtheorem{corollary}{Corollary}
\newtheorem{remark}{Remark}
\newtheorem{definition}{Definition}

\newcommand{\eqnsection}{
\renewcommand{\theequation}{\thesection.\arabic{equation}}
    \makeatletter
    \csname  @addtoreset\endcsname{equation}{section}
    \makeatother}
\eqnsection
\def\r{{\mathbb R}}
\def\e{{\mathbb E}}
\def\p{{\mathbb P}}
\def\z{{\mathbb Z}}
\def\n{{\mathbb N}}
\def\ee{\mathrm{e}}
\def\d{\, \mathrm{d}}
\def\ds{\displaystyle}
\def\eps{\varepsilon}
\def\R{{\mathbb R}}

\def\a{\alpha}
\def\b{\beta}

\def\E{{\mathbb E}}

\author[N. Enriquez]{Nathana\"el ENRIQUEZ}
\address{Laboratoire Modal'X, Universit\'e Paris 10, 200
Avenue de la R\'epublique, 92000 Nanterre, France}
\address{Laboratoire de Probabilit\'es et Mod\`eles Al\'eatoires, CNRS
UMR 7599,
Universit\'e Paris 6, 4
place Jussieu, 75252 Paris Cedex 05, France}
\email{nenriquez@u-paris10.fr}

\author[C. Sabot]{Christophe SABOT}
\address{Institut Camille Jordan, CNRS UMR 5208, Universit\'e de
Lyon, Universit\'e Lyon 1, 43, Boulevard du 11 novembre 1918,
69622 Villeurbanne Cedex} \email{sabot@math.univ-lyon1.fr}

\author[O. Zindy]{Olivier ZINDY}
\address{Laboratoire Modal'X, Universit\'e Paris 10, 200
Avenue de la R\'epublique, 92000 Nanterre, France}
\address{Laboratoire de Probabilit\'es et Mod\`eles Al\'eatoires, CNRS
UMR 7599,
Universit\'e Paris 6, 4
place Jussieu, 75252 Paris Cedex 05, France}
\email{olivier.zindy@u-paris10.fr}


\title[Limit laws for transient RWRE]{Limit laws for transient random walks in random environment on $\z$}

\begin{document}

\maketitle

\bigskip

{\footnotesize \noindent{\slshape\bfseries Abstract.} We consider
transient random walks in random environment on $\z$ with zero
asymptotic speed. A classical result of Kesten, Kozlov and Spitzer
says that the hitting time of the level $n$ converges in law, after
a proper normalization, towards a positive stable law, but they do
not obtain a description of its parameter. A different proof of this
result is presented, that leads to a complete characterization of
this stable law. The case of Dirichlet environment turns out to be
remarkably explicit.}

\bigskip
\bigskip

\section{Introduction}
   \label{s:intro}

One-dimensional random walks in random environment to the nearest
neighbors have been introduced in the sixties in order to give a
model of DNA replication. In 1975, Solomon gives, in a seminal work
\cite{solomon}, a criterion of transience-recurrence for these
walks, and shows that three different regimes can be distinguished:
the random walk may be recurrent,  or transient with a positive
asymptotic speed, but it may also be transient with zero asymptotic
speed. This  last regime, which does not exist among usual random
walks, is probably the one which is the less well understood and its
study is the purpose of the present paper.

Let us first remind the main existing results concerning the other
regimes. In his paper, Solomon computes the asymptotic speed of
transient regimes. In 1982, Sinai states, in \cite{sinai}, a limit
theorem in the recurrent case. It turns out that the motion in this
case is unusually slow since the position of the walk at time $n$
has to be normalized  by $(\log n)^2$ in order to present a non
trivial limit. In 1986, the limiting law is characterized
independently by Kesten \cite{kesten86} and Golosov \cite{golosov2}.
Let us notice here that, beyond the interest of his result, Sinai
introduces a very powerful and intuitive tool in the study of
one-dimensional random walks in random environment. This tool is the
potential, which is a function on $\z$ canonically associated to the
random environment. It turns out to be an usual random walk when the
transition probabilities at each site are independent and
identically distributed (i.i.d.).

Let us now focus on the works about the transient walk with zero
asymptotic speed. The main result was obtained by Kesten, Kozlov and
Spitzer in \cite{kesten-kozlov-spitzer} who proved  that, when
normalized by a suitable power of $n$, the hitting time of the level
$n$ converges towards a positive stable law whose index corresponds
to the power of $n$ lying in the normalization. Recently,
Mayer-Wolf, Roitershtein and Zeitouni
\cite{mayerwolf-roitershtein-zeitouni} generalized this result to
the case where the environment is defined by an irreducible Markov
chain.

Our purpose is to characterize the positive stable law in the case
of i.i.d. transition probabilities. Let us mention here that the
stable limiting law has been characterized in the case of diffusions
in random potential when the potential is either a Brownian motion
with drift \cite{kawazu-tanaka}, \cite{hu-shi-yor} or a L\'evy
process \cite{singh06a}, but we remind here that despite the
similarities of both models one cannot transport results from the
continuous model to the discrete one.

The proof chooses a radically different approach than previous ones
dealing with the transient case. The proofs in
\cite{kesten-kozlov-spitzer} and
\cite{mayerwolf-roitershtein-zeitouni} were mainly based on the
representation of the trajectory of the walk in terms of branching
processes in random environment (with immigration). 
This encoding was also used by Alili \cite{alili99b} in its study of transient 
 persistent random walks in random environment having zero asymptotic speed. 
In contrast with these works, our approach
relies heavily on Sinai's interpretation of a particle living in a
random potential. However, in the recurrent case, the potential one
has to deal with is a recurrent random walk and Sinai introduces a
notion of valley which does not make sense anymore in our setting
where the potential is a (let's say negatively) drifted random walk.
Therefore, we introduce a different notion of valley which is
closely related to the excursions of this random walk above its past
minimum. It turns out that a result of Iglehart \cite{igle}  provides the asymptotic for the distribution of the tail of the height of these excursions. Now, as
soon as one can prove that the hitting time of the level $n$ can be
reduced to the time spent by the random  walk to cross the high
excursions of the potential above its past minimum, between $0$ and
$n$, which are well separated in space, an i.i.d. property comes
out, and the problem is reduced to the study of the tail of the time
spent by the walker to cross a single excursion.

It turns out that the distribution of this tail can be expressed in terms of the expectation of the
functional of some meander associated with the random walk defining
the potential. Now, this functional is itself related to the
constant that appears in Kesten's renewal theorem \cite{kesten73}.
These last two facts are contained in \cite{enriquez-sabot-zindy+}.
Now, in the case where the transition probabilities follow some Beta
distribution a result of Chamayou and Letac \cite{chamayou-letac}
gives an explicit formula for this constant which yields finally an
explicit formula for the parameter of the positive stable law which
is obtained at the limit.

The same technics also allow to derive the convergence of the normalized process to the inverse of a standard stable subordinator. This result can be compared with the scaling limits obtained for the trap model of Bouchaud, see \cite{benarous-cerny06a} for a review.

Soon after finishing this article, we learnt of an independent work,
by Peterson and Zeitouni \cite{peterson-zeitouni}, which, by the
study of the fluctuations of the potential, showed that a quenched
stable limit law is not possible in the zero asymptotic speed
regime.

The paper is organized as follows: the results are stated in Section
\ref{s:hyp+thm},  and the rest of the paper is devoted to the
proofs.

\section{Notations and main results}
   \label{s:hyp+thm}

Let $\omega:=(\omega_i, \, i \in \z)$ be a family of i.i.d. random
variables taking values in $(0,1)$ defined on $\Omega,$ which stands
for the random environment. Denote by $P$ the distribution of
$\omega$ and by $E$ the corresponding expectation. Conditioning on
$\omega$ (i.e. choosing an environment), we define the random walk
in random environment $(X_n, \, n \ge 0)$ as a nearest-neighbor
random walk on $\z$ with transition probabilities given by $\omega$:
$(X_n, \, n \ge 0)$ is the Markov chain satisfying $X_0=0$ and for
$n \ge 0,$
$$
P_\omega \left( X_{n+1} = x+1 \, | \, X_n =x\right) = \omega_x = 1-
P_\omega \left( X_{n+1} = x-1 \, | \, X_n =x\right).
$$

\noindent We denote by $P_{\omega}$ the law of $(X_n, \, n \ge 0)$
and $E_{\omega}$ the corresponding expectation. We denote by $\p$
the joint law of $(\omega,(X_n)_{n \ge 0})$. We refer to Zeitouni
\cite{zeitouni} for an overview of results on random walks in random
environment.

In the study of one-dimensional random walks in random environment,
an important role is played by the sequence of variables 
$$
\rho_i:= \frac{1-\omega_i}{\omega_i}, \qquad i \in \z.
$$

We now introduce the hitting time $\tau(x)$ of level $x$ for the
random walk $(X_n, \, n \ge 0),$
\begin{equation}\label{hittimerw}
    \tau(x):= \inf \{ n \ge 1: \, X_n=x \}, \quad x \in \z.
\end{equation}

\noindent For $\alpha \in (0,1),$ let $\mathcal{S}_{\alpha}^{ca}$ be
a completely asymmetric (actually positive) stable random variable of index $\alpha$
with Laplace transform, for $\lambda>0,$
\begin{eqnarray*}
E[\ee^{-\lambda \mathcal{S}_{\alpha}^{ca}}]=\ee^{-
\lambda^{\alpha}}.
\end{eqnarray*}

\noindent Moreover, let us introduce the constant $C_K$ describing
the tail of Kesten's renewal series, see \cite{kesten73}, defined by
$R:=\sum_{k\geq0}\rho_0...\rho_k$:
\begin{eqnarray}
\label{tailkesten}
 P\{R>x\}\sim{C_K \over x^\kappa}, \qquad x \to
\infty.
\end{eqnarray}

\noindent Then the main result of the paper can be stated as
follows. The symbols ``$\stackrel{\mathrm{law}}{\longrightarrow}$''
 denotes the convergence in distribution.

\medskip

 \begin{theorem}
\label{t:main} Let $\omega:=(\omega_i, \, i \in \z)$ be a family of
independent and identically distributed random variables such that
\begin{itemize}
  \item[{\it (a)}] there exists $0<\kappa<1$ for which
   $E \left[  \rho_0^{\kappa} \right]=1$ and $E \left[  \rho_0^{\kappa} \log^+ \rho_0 \right]<
   \infty,$
  \item [{\it (b)}] the  distribution  of $\log \rho_0$  is non-lattice.
\end{itemize}

\noindent Then,  we have, when $n$ goes to infinity,
 \begin{eqnarray*}
  \frac{\tau(n)}{n^{1/\kappa}}  \stackrel{\mathrm{law}}{\longrightarrow} 2 \left(
{\pi\kappa^2\over\sin(\pi\kappa)}C_K^2E[\rho_0^\kappa\log\rho_0]\right)^{\frac{1}{\kappa}}
\mathcal{S}_{\kappa}^{ca},
   \\
   \\
  \frac{X_n}{n^{\kappa}}  \stackrel{\mathrm{law}}{\longrightarrow}\frac{\sin(\pi\kappa)}
{2^\kappa \pi\kappa^2 C_K^2E[\rho_0^\kappa\log\rho_0]}
\left(\frac{1}{\mathcal{S}_{\kappa}^{ca}}\right)^{\kappa}.
\end{eqnarray*}
\end{theorem}

\medskip

 \begin{remark}
Note that several probabilistic representations are available to compute $C_K$ numerically, which are equally efficient. The first one was obtained by Goldie \cite{goldie}, a second  was conjectured by Siegmund \cite{siegmund}, and we obtained a third one in \cite{enriquez-sabot-zindy+}, which plays a central role in the proof of the theorem. 
 \end{remark}

 \begin{remark}
We think that the method used in this paper could also treat the
case $\kappa=1$ (see Section \ref{sec:extension} for conjecture and
comments).
\end{remark}
\medskip

This theorem takes a remarkably explicit form in the case of Dirichlet environment, i.e. when
the law of the environment satisfies $\omega_{1}(\d
x)=\frac{1}{B(\alpha,\beta)}x^{\alpha-1}(1-x)^{\beta-1}{\bf
1}_{[0,1]}(x) \d x,$ with $\alpha, \beta >0$ and
$B(\alpha,\beta):=\int_{0}^1 x^{\alpha-1}(1-x)^{\beta-1} \d x,$
things can be made much more explicit. The assumption of Theorem
\ref{t:main} corresponds to the case where $0< \alpha-\beta <1$ and
an easy computation leads to $\kappa=\alpha-\beta.$

Thanks to a very nice result  of Chamayou and Letac \cite{chamayou-letac}
giving the explicit value of $C_K$ in this case, we obtain the following corollary:\medskip

 \begin{corollary}
\label{c:main} In the case where $\omega_1$ has a distribution
$\mathrm{Beta}(\alpha,\beta),$ with $0< \alpha-\beta <1,$ Theorem
\ref{t:main} applies with $\kappa=\alpha-\beta.$ Then, we have, when
$n$ goes to infinity,
 \begin{eqnarray*}
  \frac{\tau(n)}{n^{1/\kappa}}  \stackrel{\mathrm{law}}{\longrightarrow} 2 \left(
\frac{\pi}{\sin(\pi(\alpha-\beta))} \,
\frac{\psi(\alpha)-\psi(\beta)}{B(\alpha,\beta)^2}\right)^{\frac{1}{\alpha-\beta}}
\mathcal{S}_{\kappa}^{ca},
   \\
   \\
  \frac{X_n}{n^{\kappa}}  \stackrel{\mathrm{law}}
  {\longrightarrow}\frac{\sin(\pi(\alpha-\beta))}{2^{\alpha-\beta}\pi} \,
\frac{B(\alpha,\beta)^2}{\psi(\alpha)-\psi(\beta)}
\left(\frac{1}{\mathcal{S}_{\kappa}^{ca}}\right)^{\kappa},
\end{eqnarray*}
\noindent where $\psi$ denotes the classical Digamma function,
$\psi(z):=(\log \Gamma)'(z)=\frac{\Gamma'(z)}{\Gamma(z)}.$
\end{corollary}

\bigskip

\begin{remark} \label{r:thm3} Our technics also allow to derive the convergence of the normalized process. 
More precisely, under the assumption (a)-(b) of Theorem 1, the law of  the process
$\left(n^{-\kappa} {X_{\lfloor nt\rfloor}}, \, t \ge 0\right),$ defined on the space of c\`adl\`ag functions equipped with the uniform topology,
converges to the law of
$$
\left({\sin(\pi\kappa)\over  2^\kappa\pi \kappa^2C_K^2 E\left[\rho_0^\kappa \log\rho_0\right]}
Z_t , \, t \ge 0\right),
$$
where $Z$ is the inverse of the $\kappa$-stable subordinator $Y$ satisfying $\e[\ee^{-\lambda Y_t}] =\ee^{-t \lambda^\kappa},$ for all $\lambda>0.$ This result can be compared with the scaling limits obtained for the trap model of Bouchaud, see \cite{benarous-cerny06a} for a review.
\end{remark}
\medskip

In the following, the constant $C$ stands for a positive constant
large enough, whose value can change from line to line.

\section{Two notions of valleys}
\label{sec:prelimi} Sinai introduced in \cite{sinai} the notion of
valley in a context where the random walk defining the potential was
recurrent. We have to do a similar job in our framework where the
random walk defining the potential is negatively drifted.

Let us define precisely  the potential,
denoted by $V= (V(x), \; x\in \z)$. 
We recall first the following notation
$$
\rho_i= \frac{1-\omega_i}{\omega_i}, \qquad i \in \z.
$$

\noindent Then, the potential is a function of the environment
$\omega$ and is defined as follows:
$$
V(x) :=\left\{\begin{array}{lll} \sum_{i=1}^x \log \rho_i & {\rm if}
\ x \ge 1,
\\
 0 &  {\rm if} \ x=0,
\\
-\sum_{i=x+1}^0 \log \rho_i &{\rm if} \ x\le -1.
\end{array}
\right.
$$

\noindent Furthermore, we consider the weak descending ladder epochs
for the potential defined by $e_0:=0$ and
\begin{eqnarray*}
  e_i := \inf \{ k > e_{i-1}: \; V(k) \le V(e_{i-1})\}, \qquad i
  \ge 1,
\end{eqnarray*}

\noindent which play a crucial role in our proof. Observe that
$(e_i-e_{i-1})_{i \ge 1}$ is a family of i.i.d. random variables.
Moreover, classical results of fluctuation theory (see
\cite{feller}, p.~$396$), tell us that, under assumptions
$(a)$-$(b)$ of Theorem \ref{t:main},
\begin{eqnarray}
\label{excuinteg} E[e_1]< \infty.
\end{eqnarray}

\noindent Now, observe that the $((e_i,e_{i+1}])_{i \ge 0}$ stand
for the set of excursions of the potential above its past minimum.
Let us introduce $H_i,$ the height of the excursion $(e_i,e_{i+1}]$
defined by 
$$H_i := \max_{e_i \le k \le e_{i+1}}
\left(V(k)-V(e_i)\right),$$ 
for $i \ge 0.$ Note that the $(H_i)_{i
\ge 0}$'s are i.i.d. random variables.

The principle of the proof is to notice that the random walk in random environment spends most of its time 
climbing the high excursions. In order to quantify what  "high excursions" are, we need a key result of Iglehart  \cite{igle} which provides the asymptotic for the distribution of the tail of $H_i$, namely 
\begin{equation}\label{iglehartthm}
  \forall i\geq0,\quad   P \{H_i > h\} \sim C_I \, \ee^{- \kappa h},
     \qquad h \to \infty,
\end{equation}

\noindent where
\begin{eqnarray}\label{cst:Iglehart}
     C_I={(1-E[\ee^{\kappa V(e_1)}])^2\over \kappa
E[\rho_0^\kappa\log\rho_0]E[e_1]}.
\end{eqnarray}
Iglehart's result is actually deduced from a
former well-known result of Cramer, whose proof was later simplified
by Feller \cite{feller}, concerning the tail of the maximum  $S:=\sup \{V(k); \, k \ge 0\}$ which claims that
\begin{equation}\label{tailfeller}
P\{S>h\}\sim C_F \ee^{-\kappa h}, \qquad h \to \infty.
\end{equation}

\noindent Since $S$ is stochastically bigger than $H_0$, $C_I$ must be
smaller than $C_F$, and a rather straight argument of Iglehart shows
that the ratio between both constants is equal to $1-E[\ee^{\kappa
V(e_1)}]$.

Our strategy will be  to compute the Laplace transform of the hitting time $\tau (e_n)$ (where $\tau (x)$ is defined by (\ref{hittimerw})) which at the end will be related to $\tau (n)$ by the strong law of large numbers via $E[e_1]$.

Moreover, it appears that the times needed to cross an excursion of height $h$ is roughly of order $e^h$. Combined with Iglehart's result, it implies that the time to cross an excursion is heavy tailed for $\kappa<1$. As we know, from classical phenomena arising in the sum of heavy tailed i.i.d. random variables, the particle  will  spend most of the time at the foot  of the very few high excursions, namely those whose height has order ${\log n\over \kappa}$.
(Note that, by Iglehart's result, with an overwhelming probability, there are no excursions of height larger than ${( 1+\eps)\log n\over\kappa}$, among the $n$-first excursions.) 
This explains why the  deep valleys we define later are constructed from excursions higher than the
critical height $h_n={(1-\eps)\log n\over \kappa}$. These valleys consist actually in some portion of potential including these excursions. The high excursions are quite seldom and the valleys
are likely to be disjoint. In order to deal with almost sure
disjoint valleys, we also introduce $*$-valleys which coincide with
deep valleys with high probability.

\subsection{The deep valleys}
\label{deepvalleys}

\noindent Let us define the maximal variations of the potential
before site $x$ by
\begin{eqnarray*}
 \label{incrempoten}
    V^{\uparrow}(x):= \max_{0 \le i \le j \le x} (V(j)-V(i)), \qquad x \in
    \n,
    \\
    V^{\downarrow}(x):= \min_{0 \le i \le j \le x} (V(j)-V(i)), \qquad x \in
    \n.
\end{eqnarray*}
By extension, we introduce
\begin{eqnarray*}
 \label{incremupgene}
    V^{\uparrow}(x,y):= \max_{x \le i \le j \le y} (V(j)-V(i)),
    \qquad  x<y,
    \\
\label{incremdowngene}
    V^{\downarrow}(x,y):= \min_{x \le i \le j \le y} (V(j)-V(i)),
    \qquad x<y.
\end{eqnarray*}

In order to define deep valleys, we  extract from the first $n$
excursions of the potential above its minimum, these whose heights
are greater than a critical height $h_n,$ defined by
\begin{equation}
\label{hcritic} h_n:= {(1-\varepsilon) \over \kappa} \log n,
\end{equation}
for some $0<\varepsilon<1/3,$ see Figure 1.   Let $(\sigma(i))_{i \ge 1}$
be the successive indexes of excursions, whose heights are greater
than $h_n.$ More precisely,
\begin{eqnarray*}
\sigma(1)&:=&\inf \{ i \ge 0 : H_i \ge h_n,\},
\\
\sigma(j)&:=& \inf \{ i >\sigma(j-1): H_i \ge h_n \}, \qquad j \ge
2,
\\
K_n&:=& \max \{ j\ge 0 :\sigma(j)\le n \}.
\end{eqnarray*}

   \begin{figure}[htb]
    \label{fig1}
\begin{tikzpicture}
\draw[thick][->] (0.5,0)--(14.5,0); 
\draw[thick][->] (1.5,-3.8)--(1.5,1);
\draw[thick](1.5,0)--(1.0,0.16);

\draw[thick](1.5,0)--(4.5,-1);
\draw[thick](4.5,-1)--(4.7,-0.85);
\draw[thick](4.7,-0.85)--(6.7,-1.5);
\draw[thick](6.7,-1.5)--(6.95,-1.32);
\draw[thick](6.7,-1.5)--(6.95,-1.32);
\draw[thick](6.95,-1.32)--(8.45,-1.81);
\draw[thick](8.65,-1.87)--(8.85,-1.93);
\draw[thick](9.05,-1.99)--(9.25,-2.05);
\draw[thick](9.45,-2.11)--(11.45, -2.77);
\draw[thick](11.45, -2.77)--(11.65, -2.62);
\draw[thick](11.65, -2.62)--(14.2,-3.5);

\draw (4.4,-1.25) node {$b_1$}; 
\draw (4.8,-0.65) node {$c_1$}; 
\draw (6.6,-1.725) node {$b_2$}; 
\draw(7.05,-1.1) node {$c_2$}; 
\draw (11.4,-3.015) node {$b_{K_n}$}; 
\draw(11.8,-2.39) node {$c_{K_n}$}; 

%\draw (12,-1) node {$N=N(n)$};

\draw (14.4,0.3) node {$x$};
\draw (0.9,0.9) node {$V(x)$}; 
\draw(14.1,-0.25) node {$e_n$}; 
\draw[thick](14.1,-0.1)--(14.1,0.1);
\draw (1.3,-0.25) node {$0$};

% \draw[thick](6,-1.05)--(6,-0.95);
 %\draw[thick](6.25,-0.85)--(6.25,-0.75);

\end{tikzpicture}
\caption{Potential and valleys.}
\end{figure}
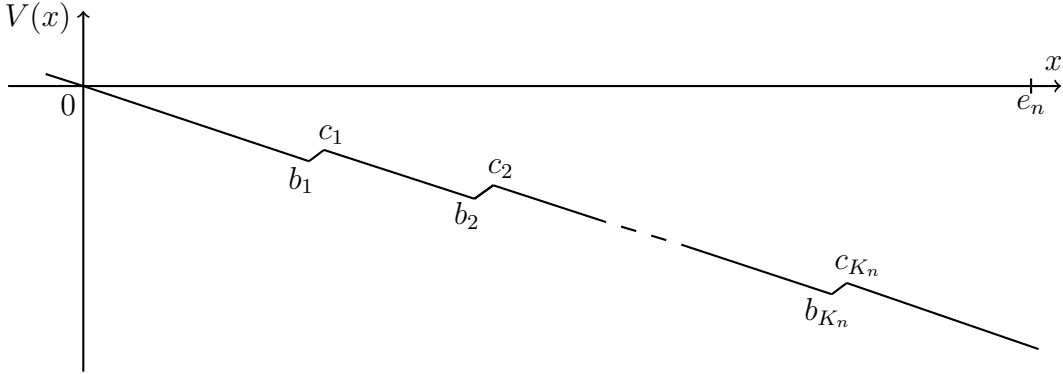

\noindent We consider now some random variables depending only on $n$ and on
the environment, which define the deep valleys.

\medskip

\begin{definition}
 \label{defvalley}
For $1\le j \le K_n+1,$ let us introduce
\begin{eqnarray*}
b_j&:=&e_{\sigma(j)},
\\
a_j&:=&\sup \{ k \le b_j : \, V(k)-V(b_j) \ge D_n  \},
\\
T_{j}^{\uparrow}&:=&\inf \{ k \ge b_j : \, V(k)-V(b_j) \ge h_n  \},
\\
\overline{d}_j&:=&e_{\sigma(j)+1},
\\
c_j&:=&\inf \{ k \ge b_j : \, V(k)=\max_{b_j \le x \le
\overline{d}_j} V(x) \},
\\
d_j&:=&\inf \{ k \ge \overline{d}_j: \, V(k)-V(\overline{d}_j) \le
-D_n \}.
\end{eqnarray*}
where $D_n:=(1+{1 \over \kappa}) \log n.$ We call
$(a_j,b_j,c_j,d_j)$ a deep valley and denote by $H^{(j)}$ the height
of the $j$-th deep valley.
\end{definition}
Note that all the random variables introduced in this section depend on $n,$ see Figure 2.
\medskip

\begin{remark}
 \label{2valleys}
It may happen that two different deep valleys are not disjoint, even
if this event is highly improbable as it will be shown in Lemma
\ref{l:Aboxgap} and Lemma \ref{l:Aboxlength} in Subsection
\ref{themule}.
\end{remark}

   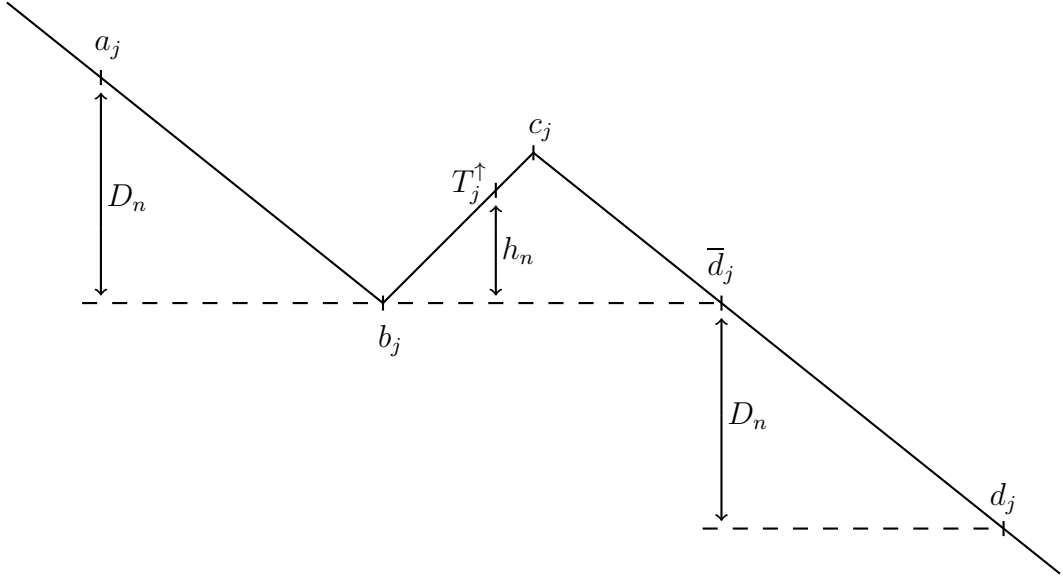
\begin{figure}[htb]
\begin{tikzpicture}
   \label{fig2}

\draw[thick](1.5,-1)--(6.5,-5);
\draw[thick](6.5,-5)--(8.5,-3);
\draw[thick](8.5,-3)--(15.5,-8.6);

\draw[thick](6.3,-5)--(6.1,-5);
\draw[thick](5.7,-5)--(5.9,-5);
\draw[thick](5.5,-5)--(5.3,-5);
\draw[thick](4.9,-5)--(5.1,-5);
\draw[thick](4.7,-5)--(4.5,-5);
\draw[thick](4.1,-5)--(4.3,-5);
\draw[thick](3.7,-5)--(3.9,-5);
\draw[thick](3.3,-5)--(3.5,-5);
\draw[thick](3.1,-5)--(2.9,-5);
\draw[thick](2.5,-5)--(2.7,-5);

\draw[thick](6.7,-5)--(6.9,-5);
\draw[thick](7.1,-5)--(7.3,-5);
\draw[thick](7.5,-5)--(7.7,-5);
\draw[thick](7.9,-5)--(8.1,-5);
\draw[thick](8.3,-5)--(8.5,-5);
\draw[thick](8.7,-5)--(8.9,-5);
\draw[thick](9.1,-5)--(9.3,-5);
\draw[thick](9.5,-5)--(9.7,-5);
\draw[thick](9.9,-5)--(10.1,-5);
\draw[thick](10.3,-5)--(10.5,-5);
\draw[thick](10.7,-5)--(10.9,-5);

\draw[thick](6.5,-5.1)--(6.5,-4.9);
\draw[thick](8.5,-3.1)--(8.5,-2.9);
\draw[thick](8,-3.6)--(8,-3.4);
\draw[thick](11,-5.1)--(11,-4.9);
\draw[thick](2.75,-2.1)--(2.75,-1.9);
\draw[thick](14.75,-8.1)--(14.75,-7.9);

\draw[thick](14.35,-8)--(14.55,-8);
\draw[thick](14.15,-8)--(13.95,-8);
\draw[thick](13.75,-8)--(13.55,-8);
\draw[thick](13.35,-8)--(13.15,-8);
\draw[thick](12.95,-8)--(12.75,-8);
\draw[thick](12.55,-8)--(12.35,-8);
\draw[thick](12.15,-8)--(11.95,-8);
\draw[thick](11.75,-8)--(11.55,-8);
\draw[thick](11.35,-8)--(11.15,-8);
\draw[thick](10.95,-8)--(10.75,-8);

\draw[thick][->] (8,-4)--(8,-3.7); 
\draw[thick][->] (8,-4)--(8,-4.9);

\draw[thick][->] (2.75,-4)--(2.75,-2.2); 
\draw[thick][->] (2.75,-4)--(2.75,-4.9);

\draw[thick][->] (11,-6.5)--(11,-7.9); 
\draw[thick][->] (11,-6.5)--(11,-5.2);

\draw (6.6,-5.5) node {$b_j$}; 
\draw (8.6,-2.7) node {$c_j$}; 
\draw (2.85,-1.6) node {$a_j$}; 
\draw(14.75,-7.6)node {$d_j$}; 
\draw (11,-4.5) node {$\overline d _j$}; 
\draw(7.65,-3.4) node {$T^{\uparrow}_j$};

\draw (8.3,-4.3) node {$h_n$}; 
\draw (11.35,-6.5) node {$D_n$}; 
\draw (3.1,-3.6) node {$D_n$}; 

\end{tikzpicture}
\caption{Zoom on the $j$-th valley.}
\end{figure}

\subsection{The $*$-valleys}
\label{*-valleys}

Let us introduce now a subsequence of the deep valleys defined
above. It will turn out that both sequences coincide with
probability tending to $1$ as $n$ goes to infinity. This  will be
specified in Lemma \ref{model1+2}. Let us first introduce
\begin{eqnarray*}
\gamma^*_1&:=&\inf \{ k \ge 0 : \, V(k) \le -D_n  \},
\\
T_{1}^{*}&:=&\inf \{ k \ge \gamma^*_1 : \,
V^{\uparrow}(\gamma^*_1,k) \ge h_n \},
\\
b^*_1&:=&\sup \{ k \le T_{1}^{*} : \, V(k)=\min_{0 \le x \le
T_{1}^{*}} V(x) \},
\\
a^*_1&:=&\sup \{ k \le b^*_1 : \, V(k)-V(b^*_1) \ge D_n  \},
\\
\overline{d}^*_1&:=&\inf \{ k \ge T_{1}^{*} : \, V(k) \le
V(b_1^*)\},
\\
c^*_1&:=&\inf \{ k \ge b^*_1 : \, V(k)=\max_{b^*_1\le x \le
\overline{d}^*_1} V(x) \},
\\
d^*_1&:=&\inf \{ k \ge \overline{d}^*_1 : \,
V(k)-V(\overline{d}^*_1) \le -D_n  \}.
\end{eqnarray*}

\noindent Let us define the following sextuplets of points by
iteration
$$
(\gamma^*_j,a^*_j,b^*_j,T_{j}^{*},c^*_j,\overline{d}^*_j,d^*_j):=
(\gamma^*_1,a^*_1,b^*_1,T_{1}^{*},c^*_1,\overline{d}^*_1,d^*_1)
\circ \theta_{d_{j-1}^*}, \qquad j \ge 2,
$$
where $\theta_{i}$ denotes the $i$-shift operator.

\begin{definition}
 \label{def*valley}
We call a $*$-valley any quadruplet $(a^*_j,b^*_j,c^*_j,d^*_j)$ for
$j \ge 1.$ Moreover, we shall denote by $K_n^*$ the number of such
$*$-valleys before $e_n,$ i.e. $ K_n^*:= \sup \{ j \ge 0 : T_{j}^{*}
\le e_n \}.$
\end{definition}

\noindent It will be made of independent and identically distributed
portions of potential (up to some translation).

\section{Reduction to a single valley}
\label{sec:part1}

This section is devoted to the proof of Proposition \ref{P:prop1}
which tells that the study of $\tau(e_n)$ can be reduced to the
analysis of the time spent by the random walk to cross the first
deep valley. To ease notations, we introduce $\lambda_n :={\lambda
\over n^{1/\kappa}}.$

\medskip

\begin{proposition}
 \label{P:prop1} For all $n$ large enough, we have
 \begin{eqnarray*}
  \e \left[ \ee^{- \lambda_n \, \tau(e_n)} \right] &\in&
  \left[ E \left[ E_{\omega,|a_1}^{b_1} \left[ \ee^{- \lambda_n
  \tau(d_1)} \right] \right]^{\overline{K}_n} +o(1) \, , \,   E \left[
  E_{\omega,|a_1}^{b_1} \left[ \ee^{- \lambda_n \tau(d_1)} \right]
  \right]^{\underline{K}_n}+o(1) \right].
\end{eqnarray*}

\noindent where $\underline{K}_n:= \lfloor n q_n
(1-n^{-\varepsilon/4}) \rfloor,$ $\overline{K}_n:= \lceil n q_n
(1+n^{-\varepsilon/4}) \rceil,$ $q_n:=P\{H_0\ge h_n\}$ and where
$E_{\omega,|y}^x$ denotes the quenched law of the random walk in the
environment $\omega,$ starting at $x$ and reflected at site $y$.
\end{proposition}

\medskip

\subsection{Introducing ``good'' environments}
\label{themule}

Let us define the four following events, that concern exclusively
the potential $V.$ The purpose of this subsection is to show that
they are realized with an asymptotically overwhelming probability
when $n$ goes to infinity. These results will then make it possible
to restrict the study of $\tau(e_n)$ to these events.
\begin{eqnarray*}
  A_1(n) &:=& \left\{ e_n <C' n \right\},
  \\
  A_2(n) &:=& \left\{\lfloor n q_n(1-n^{-\varepsilon/4}) \rfloor \le
  K_n \le \lceil n q_n(1+n^{-\varepsilon/4}) \rceil \right\},
     \\
  A_3(n) &:=& \cap_{j=0}^{K_n} \left\{ \sigma(j+1)-\sigma(j) \ge n^{1-3\varepsilon}
    \right\},
     \\
  A_4(n) &:=& \cap_{j=1}^{K_n+1} \left\{ d_j- a_j
    \le C'' \log n \right\},
\end{eqnarray*}

\noindent where $\sigma(0):=0$ (for convenience of notation) and
$C',$ $C''$ stand for positive constants which will be specified
below.

In words, $A_1(n)$ allows us to bound the total length of the first
$n$ excursions. The event $A_2(n)$ gives a control on the number of
deep valleys. The event $A_3(n)$ ensures that the deep valleys are
well separated, while $A_4(n)$ bounds finely the length of each of
them. Before proving that the $A_i$'s are typical events, 
let us first give a preliminary result concerning large deviations that we will use throughout the paper.\medskip

\begin{lemma}
 \label{lemmeLD}
 Under assumption $(a)$, large deviations occur for the potential seen as a sum of i.i.d. random variables. Indeed for all $x\ge m:=E \left[\log \rho_0 \right]$ (recall that $(a)$ implies $m<0$) and all $j\ge 1,$ we have
 \begin{eqnarray}
  \label{eq:lemmeLD1}
P \left\{V(j) \ge j x\right\} \le \exp\{-j I(x)\},
\end{eqnarray}
 with $I(x):=\sup_{t \ge 0}\{tx-\Lambda(t)\}$ and $\Lambda(t):=\log E[\rho_0^t].$ Moreover, the rate function $I$ is lower semicontinuous, satisfies $I(0)>0$ and 
\begin{eqnarray}
    \inf_{x > 0} {I(x) \over x} \ge \kappa.
     \label{S3+}
\end{eqnarray}
\end{lemma}

\medskip

\begin{proof} 
Let us first prove (\ref{eq:lemmeLD1}) which is the upper bound in Cramer's theorem in $\r,$ see \cite{dembo-zeitouni98}. Observe first that for all $x$ and every $t \ge 0,$ an application of Markov's inequality yields
\begin{eqnarray}
P\{V(j)  \ge j x\} &=&  E [{\bf 1}_{\{ V(j)-jx\ge0 \}}] \le E [\ee^{t(V(j)-jx)}]
\nonumber
\\
&=& \ee^{-j t x} E [\ee^{t \log \rho_0}]^j=\ee^{-j  \{ tx - \Lambda(t) \}}.
\label{eq:ld12}
\end{eqnarray}
 Then, we get (\ref{eq:lemmeLD1})  by taking the infimum over $t\ge0$ in (\ref{eq:ld12}).

To prove that $I(0)>0,$ observe first that $I(0)=- \inf_{t \ge 0} \Lambda(t).$
Now since the function $g(t):=E[\rho_0^t]$ satisfies $g(0)=g(\kappa)=1$ (by assumption $(a)$)  and $g'(0)<0$ (indeed $g'(0)=E \left[\log \rho_0 \right]<0$), we get that $\inf_{0 \le t \le \kappa} g(t) <1,$ which implies  $-\inf_{t \ge 0} \Lambda(t)>0.$

The proof of (\ref{S3+})  is straightforward.
Indeed, recalling that
$I(x)=\sup_{t \ge 0}\{tx-\Lambda(t)\}$ for $x > 0,$ we have
$I(x)\ge \kappa x -\Lambda(\kappa)=\kappa x,$ since
$\Lambda(\kappa)=0.$  \end{proof}
\medskip
Note that the claim of (\ref{S3+}) appears on page 236 in \cite{zeitouni} and that \cite{zeitouni} claims an equality under certain assumptions.

 Now, let us introduce the following hitting times (for the
potential)
\begin{eqnarray*}
\label{hittimepot1}
    T_{h}:= \min \{x \ge 0: \, V(x) \ge h \}, \qquad
    h>0,
    \\
\label{hittimepot2}
    T_{A}:= \min \{x \ge 0: \, V(x) \in A \}, \qquad
    A \subset \r.
\end{eqnarray*}
and prove that the $A_i(n)$'s occur with an overwhelming probability when $n$ tends to infinity.

\medskip
\begin{lemma}
 \label{l:totallenght}
The probability $P\{ A_1(n)\}$ converges to $1$ when $n$ goes to
infinity.
\end{lemma}

\begin{proof} It is a direct consequence of the law of
large numbers as soon as $C'$ is taken bigger than $E[e_1].$
\end{proof}

\medskip

\begin{lemma}
 \label{l:Aboxnumber}
The probability $P\{ A_2(n)\}$ converges to $1$ when $n$ goes to
infinity.
\end{lemma}

\noindent  In words, Lemma \ref{l:Aboxnumber} means that $K_n$
``behaves" like $C_I n^{\varepsilon},$ when $n$ tends to infinity.
In particular, (\ref{iglehartthm}), which yields $q_n \sim {C_I
\over n^{1-\varepsilon}},$ and Lemma \ref{l:Aboxnumber} imply
\begin{eqnarray}
\label{Knmajo}
  P\{ K_n+1 \ge 2 C_I n^{\varepsilon}\} \to 0, \qquad n \to \infty.
\end{eqnarray}

\begin{proof} At first, observe that
$$P\big\{ {K_n \over n  q_n} \ge 1+n^{-\varepsilon/4}\big\}
=P\{ K_n - n  q_n \ge n^{1-\varepsilon/4} q_n \} \le {{\rm Var}(K_n)
\over n^{2(1-\varepsilon/4)} q_n^2},
$$

\noindent the inequality being a consequence of Markov inequality
and the fact that $K_n$ follows a binomial distribution of parameter
$(n,q_n).$ Moreover, ${\rm Var}(K_n)=n q_n (1-q_n) \le n q_n$
implies
$$P\{ {K_n \over n  q_n} \ge 1+n^{-\varepsilon/4}\}
 \le {1 \over n^{1-\varepsilon/2} q_n}.
$$

\noindent Now, Iglehart's result (see (\ref{iglehartthm})) implies
$q_n \sim {C_I \over n^{1-\varepsilon}},$ $n \to \infty.$ Therefore
we get that $P\{ {K_n \over n  q_n} \le 1+n^{-\varepsilon/4} \}$
converges to $1$ when $n$ goes to infinity. Using similar arguments,
we get the convergence to $1$ of $P\{ {K_n \over n q_n} \ge
1-n^{-\varepsilon/4}\}.$  \end{proof}

\medskip
\begin{lemma}
 \label{l:Aboxgap} The probability $P\{ A_3(n)\}$ converges to $1$ when $n$
goes to infinity.
\end{lemma}

\begin{proof} We make first the trivial observation that
\begin{eqnarray*}
P\{A_3(n)\} &\ge& P\{\sigma(j+1)-\sigma(j) \ge n^{1-3\varepsilon},
\,0 \le j \le \lfloor 2 C_I n^{\varepsilon}\rfloor
\, ; \, K_n \le 2 C_I n^{\varepsilon}\} \\
&\ge& P\{\sigma(j+1)-\sigma(j) \ge n^{1-3\varepsilon}, \, 0 \le j
\le \lfloor2 C_I n^{\varepsilon}\rfloor\} - P\{K_n \ge 2 C_I
n^{\varepsilon}\},
\end{eqnarray*}

\noindent the second inequality being a consequence of $P\{A \, ;
B\}\ge P\{A\}-P\{B^c\},$ for any couple of events $A$ and $B.$
Therefore, recalling (\ref{Knmajo}) and using the fact that
$(\sigma(j+1)-\sigma(j))_{0 \le j \le \lfloor 2 C_I
n^{\varepsilon}\rfloor}$ are i.i.d. random variables, it remains to
prove that
\begin{eqnarray*}
\label{Aboxgapeq1}
    P\{\sigma(1) \ge
n^{1-3\varepsilon}\}^{\lfloor2 C_I n^{\varepsilon}\rfloor} \to 1,
\qquad n \to \infty.
\end{eqnarray*}

\noindent Since $\sigma(1)$ is a geometrical random variable with
parameter $q_n,$ $P\{\sigma(1) \ge n^{1-3\varepsilon}\}$ is equal to
$(1-q_n)^{\lceil n^{1-3\varepsilon} \rceil},$ which implies
\begin{eqnarray*}
     P\{\sigma(1) \ge
n^{1-3\varepsilon}\}^{\lfloor2 C_I
n^{\varepsilon}\rfloor}=(1-q_n)^{\lfloor2 C_I n^{\varepsilon}\rfloor
\, \lceil n^{1-3\varepsilon} \rceil} \ge \exp \left\{-C
n^{1-2\varepsilon}  q_n \right\}.
\end{eqnarray*}

\noindent Then, the conclusion follows from (\ref{iglehartthm}),
which implies that $q_n \sim C_I / n^{1-\varepsilon},$ $n \to
\infty.$ \end{proof}

\medskip
\begin{lemma}
 \label{l:Aboxlength} For $C''$ large enough, the probability $P\{ A_4(n)\}$
converges to $1$ when $n$ goes to infinity.
\end{lemma}

\begin{proof} Looking at the proof of Lemma
\ref{l:Aboxgap}, we have to prove that $P\{d_j- a_j \ge C'' \log n
\}$ is equal to a $o(n^{-\varepsilon}),$ $n \to \infty.$ Moreover,
observing that $d_j- a_j=(d_j-\overline{d}_j)+(\overline{d}_j-
T_{j}^{\uparrow}) +(T_{j}^{\uparrow}-b_j)+ (b_j - a_j),$ the proof
of Lemma \ref{l:Aboxlength} boils down to showing that, for $C''$
large enough,
\begin{eqnarray}
  P\{d_j-\overline{d}_j \ge \frac{C''}{4} \log n\}
  &=& o(n^{-\varepsilon}), \qquad n \to
  \infty,
  \label{Aboxlengtheq0}
  \\
  P\{\overline{d}_j-T_{j}^{\uparrow} \ge \frac{C''}{4} \log n \}
  &=& o(n^{-\varepsilon}), \qquad n \to
  \infty,
  \label{Aboxlengtheq1}
  \\
  P\{T_{j}^{\uparrow}-b_j \ge \frac{C''}{4} \log n\}
  &=&o(n^{-\varepsilon}), \qquad n \to
  \infty,
  \label{Aboxlengtheq2}
  \\
  P\{b_j - a_j \ge \frac{C''}{4} \log n\}
  &=&o(n^{-\varepsilon}),
\qquad n \to \infty.
  \label{Aboxlengtheq3}
\end{eqnarray}

To prove (\ref{Aboxlengtheq0}), we apply the strong Markov property
at time $\overline{d}_j$ such that we get $P\{d_j-\overline{d}_j \ge
\frac{C''}{4} \log n\} \le P\{T_{(-\infty,-D_n]} \ge \frac{C''}{4}
\log n\}.$ Therefore, we have
\begin{eqnarray*}
P\{d_j-\overline{d}_j \ge \frac{C''}{4} \log n\} \le P\{ \inf_{0 \le
x \le \frac{C''}{4} \log n} V(x) > -D_n \} \le P\{ V(\frac{C''}{4}
\log n) > -D_n\}.
\end{eqnarray*}

\noindent Recalling that $D_n:=(1+{1 \over \kappa}) \log n,$ we can use Lemma \ref{lemmeLD}, which implies $P\{ V(\frac{C''}{4} \log n)
> -D_n\}\le \ee^{ - \frac{C''}{4} \log n \, I(- {4 \over C''}(1+ {1 \over
\kappa}))}.$ Then, this inequality implies (\ref{Aboxlengtheq0}) by choosing $C''$
large enough such that $\frac{C''}{4} \, I(- {4 \over C''}(1+ {1
\over \kappa}))> \varepsilon,$ which is possible since $I(0)>0.$

To prove (\ref{Aboxlengtheq1}), observe first that
(\ref{iglehartthm}) implies $ P\{ H^{(j)}> {(1+\varepsilon') \over
\kappa} \log n \} \sim
n^{-(\varepsilon'+\varepsilon)}=o(n^{-\varepsilon}),$ $n \to
\infty.$ Therefore, we obtain that
$P\{\overline{d}_j-T_{j}^{\uparrow} \ge \frac{C''}{4} \log n\}$ is
less or equal than $P\{T_{(-\infty,-{1+\varepsilon' \over \kappa}
\log n]} \ge \frac{C''}{4} \log n \}  + o(n^{-\varepsilon})$ and
conclude the proof with the same arguments we used to treat
(\ref{Aboxlengtheq0}).

To get (\ref{Aboxlengtheq2}), observe first that
\begin{eqnarray*}
  P\{T_{j}^{\uparrow}-b_j \ge \frac{C''}{4} \log n
\}&=&P\{T_{h_n} \ge \frac{C''}{4} \log n \, | \, H_0 \ge h_n \}
  \\
   &\le& P\{ \frac{C''}{4} \log n \le T_{h_n}< \infty
\}/P\{H_0 \ge h_n\}.
\end{eqnarray*}

\noindent Furthermore, Lemma \ref{lemmeLD} 
yields
\begin{eqnarray*}
  P\{ \frac{C''}{4} \log n \le T_{h_n}< \infty
\} &\le& \sum_{k \ge \frac{C''}{4} \log n}  P\{ V(k) \ge h_n \}
 \le \sum_{k \ge \frac{C''}{4} \log n}  \ee^{-k \, I\left({h_n \over k}\right)}
 \\
   &\le & \sum_{k \ge \frac{C''}{4} \log n}  \ee^{-k \, I(0)} \le {C \over n^{\frac{C''}{4} \,
   I(0)}},
\end{eqnarray*}

\noindent the third inequality being a consequence of the fact that
the convex rate function $I(\cdot)$ is an increasing function on
$(m,+\infty).$ Using (\ref{iglehartthm}), we get, for all large $n,$
$$
P\{T_{j}^{\uparrow}-b_j \ge \frac{C''}{4} \log n\} \le {C \over
n^{\frac{C''}{4} I(0)- (1-\varepsilon)}},
$$

\noindent which yields (\ref{Aboxlengtheq2}), by choosing $C''$
large enough such that $C''>{4 \over I(0)}.$

For (\ref{Aboxlengtheq3}), observe first that
$((V(k-b_j)-V(b_j))_{a_j \le k \le b_j},a_j,b_j)$ has the same
distribution as $((V(k))_{a^- \le k \le 0},a^-,0)$ under $P\{\cdot
|V(k)\ge0, a^- \le k \le 0 \},$ where $a^-:=\sup \{ k \le 0 : \,
V(k) \ge D_n  \}$. Then, since $P\{V(k)\ge0, k \le 0 \}>0$ and since
$(V(-k), \, k \ge 0)$ has the same distribution as $(-V(k), \, k \ge
0),$ we obtain
$$
P\{b_j - a_j \ge \frac{C''}{4} \log n\} \le C P\{T_{(-\infty,-D_n]}
> \frac{C''}{4} \log n\} \le C P\{ V(\frac{C''}{4} \log n) > -D_n\}.
$$

\noindent Now, the arguments are the same as in the proof of
(\ref{Aboxlengtheq0}). \end{proof}

\bigskip

Defining  $A(n):=A_1(n) \cap A_2(n) \cap A_3(n) \cap A_4(n),$ a
consequence of Lemma \ref{l:totallenght}, Lemma \ref{l:Aboxnumber},
Lemma \ref{l:Aboxgap} and Lemma \ref{l:Aboxlength}, is that
\begin{eqnarray}
\label{P:mule1}
  P\{ A(n)\} \to 1.
\end{eqnarray}

\noindent The following lemma tells us that the $*$-valleys coincide
with the sequence of deep valleys with an overwhelming probability
when $n$ goes to infinity.

\medskip

\begin{lemma}
\label{model1+2} If $A^*(n):=\{  K_n=K_n^* \, ; \,
(a_j,b_j,c_j,d_j)= (a^*_j,b^*_j,c^*_j,d^*_j), \, 1 \le j \le K_n
\},$ then we have that the probability $P\{A^*(n)\}$ converges to
$1,$ when $n$ goes to infinity.
\end{lemma}

\begin{proof} Since, by definition, the $*$-valleys
constitute a subsequence of the deep valleys, Lemma \ref{model1+2}
is a consequence of Lemma \ref{l:Aboxgap} together with Lemma
\ref{l:Aboxlength}. \end{proof}

\medskip

\begin{remark}
Another meaning of this result is that, with probability tending to
$1,$ two deep valleys are necessarily disjoint.
\end{remark}

\medskip

\subsection{Preparatory lemmas}

In this subsection, we develop some technical tools allowing us to
improve our understanding of the random walk's behavior. In Lemma
\ref{l:ADT}, we prove that, after exiting a deep valley, the random
walk will not come back to another deep valley it has already
visited, with probability tending to one. Moreover, Lemma
\ref{l:ADTgamma} specifies that the random walk typically exits from
a $*$-valley on the right, while Lemma \ref{l:IA} shows that the
time spent between two deep valleys is negligible. Lemma \ref{l:valley-*valley} 
 states that the first valley coincides with the first $*$-valley with probability $1-o(n^{-\varepsilon}),$ when $n$ goes to infinity.

\subsubsection{Preliminary estimates for inter-arrival times}
\label{IAtimes}

Let us introduce
\begin{eqnarray*}
\label{firstincrempot}
    T^{\uparrow}(h):= \min \{x \ge 0: \, V^{\uparrow}(x) \ge h \}, \qquad
    h >0,
    \\
    T^{\downarrow}(h):= \min \{x \ge 0: \, V^{\downarrow}(x) \le -h \}, \qquad
    h >0.
\end{eqnarray*}

\medskip

\begin{lemma}
 \label{lemmaenv}
Under assumptions of Theorem \ref{t:main}, we have, for $h$ large enough,
\begin{eqnarray*}
\label{majoenv}
    \e_{|0}\left[\tau_h \right] \le C \,  \ee^{h},
\end{eqnarray*}
\noindent where $\e_{|0}$ denotes the expectation under the law
$\p_{|0}$ of the random walk in the random environment $\omega$
(under $P$) reflected at $0$ and $\tau_h:=\tau(T^{\uparrow}(h)-1).$
\end{lemma}

\begin{proof} Using (Zeitouni \cite{zeitouni}, formula
(2.1.14)), we obtain that $ \e_{|0} \left[\tau_h\right]$ is bounded
from above by $E\big[ \sum_{0 \le i \le j < T^{\uparrow}(h)}
\ee^{V(j)-V(i)} \big].$ Therefore, since $T^{\uparrow}(h) \le
T^{\uparrow}(h) \circ \theta_i$ for any $i \ge 0$ (where $\theta$ denotes the shift operator for the environment), we obtain
\begin{eqnarray}
\e_{|0} \left[\tau_h\right] \le \sum_{i \ge 0} E\Big[ {\bf
1}_{\{i<T^{\uparrow}(h)\}} \sum_{i \le j < T^{\uparrow}(h)}
\ee^{V(j)-V(i)} \Big] \le \beta_1(h) \, \beta_2(h), \label{expo:eq1}
\end{eqnarray}
where
\begin{eqnarray*}
\beta_1(h) &:=& E\big[T^{\uparrow}(h)\big],
\\
\beta_2(h) &:=& E\Big[ \sum_{0 \le j < T^{\uparrow}(h)} \ee^{V(j)}
\Big].
\end{eqnarray*}
To bound $\beta_1(h),$ let us introduce the number $N$ of complete
excursions before $T^{\uparrow}(h),$ defined by $N=N(h):=\sup\{i \ge
0: e_i<T^{\uparrow}(h) \}.$ Then, we can write
$\beta_1(h)=E[\sum_{i=0}^{N-1}(e_{i}-e_{i-1})+(T^{\uparrow}(h)-e_N)].$
Observe that the definition of $T^{\uparrow}(h)$ implies that $N$ is
a geometrical random variable with parameter $q=q(h):=P\{H \ge h\}$
and recall that, by (\ref{iglehartthm}), we have $q \sim C_I \,
\ee^{- \kappa h},$ $h \to \infty.$  Therefore, we get, for $h$ large
enough,
\begin{eqnarray*}
\beta_1(h) &\le& \sum_{k \ge 0} (1-q)^k q \big(k E[e_1 | H <
h]+E[T_h | H \ge h]\big)
 \\
&\le& C \sum_{k \ge 0} (1-q)^k q \big(k E[e_1]+ E[T_h | H \ge
h]\big),
\end{eqnarray*}
the second inequality being a consequence of the fact that
$E[e_1]<\infty$ (see (\ref{excuinteg})) together with $P\{H< h\} \to
1,$ $h \to \infty,$ by (\ref{iglehartthm}). By obvious calculations,
this yields $\beta_1(h) \le C (1-q) q^{-1} E[e_1]+E[T_h | H \ge h],$
which implies with (\ref{iglehartthm}) that
\begin{eqnarray}
\label{expo:eq2} \beta_1(h) \le C \ee^{\kappa h}+E[T_h | H \ge h].
\end{eqnarray}
Now, let us bound $E[T_h | H \ge h].$ For this purpose, we observe first
that $E[T_h | H \ge h] \le C \ee^{\kappa h} \sum_{k \ge 0} (k+1)
P\{T_h=k+1 \,;\, H \ge h \}.$ Then, applying the Markov property at
time $k,$ we get
\begin{eqnarray*}
E[T_h | H \ge h] &\le& C \ee^{\kappa h} \sum_{k \ge 0} (k+1) E[{\bf
1}_{\{0<V(k)<h\}} \ee^{-\kappa(h-V(k))}]
\\
&\le& C\sum_{k \ge 0} (k+1)\sum_{j=0}^{\lfloor h \rfloor}
\ee^{\kappa (j+1)} P\{V(k)\ge j\}.
\end{eqnarray*}
By Lemma \ref{lemmeLD}, we have $P\{V(k)\ge j\}
\le \ee^{-k I(\frac{j}{k})}.$ Now, the fact that $I(\cdot)$ is an
increasing function on $\r^+$ along with (\ref{S3+}) imply
$$
P\{V(k)\ge j\} \le  \ee^{-{k\over 2} I(\frac{j}{k})}  \ee^{-{k\over 2} I(\frac{j}{k})}  \le \ee^{- k \frac{I(0)}{2}} \ee^{-\kappa
\frac{j}{2}}.
$$
Since $I(0)>0$, this yields that there exists $C>0$ such that, for
all large $h,$
\begin{eqnarray}
\label{expo:eq3} E[T_h | H \ge h] \le C \ee^{\frac{\kappa}{2} h}.
\end{eqnarray}
Combining together (\ref{expo:eq2}) and (\ref{expo:eq3}), we obtain for $h$
large enough,
\begin{eqnarray}
\label{expo:eq4}
 \beta_1(h) \le C \ee^{\kappa h}.
\end{eqnarray}

Let us now bound $\beta_2(h).$ We introduce first
$\mathcal{E}_k:=\{ \max_{0 \le j \le k-1}H_j <h   \,;\, H_k \ge h
\}$ and write
\begin{eqnarray*}
 \beta_2(h)&=&\sum_{k \ge 0} E \Big[{\bf
1}_{\mathcal{E}_k} \sum_{0 \le j < T^{\uparrow}(h)} \ee^{V(j)}\Big]
\\
&=& \sum_{k \ge 0}\Big(\sum_{i=0}^{k-1} E\Big[{\bf
1}_{\mathcal{E}_k}
 \ee^{V(e_j)} J_{i}\Big] + E\Big[{\bf
1}_{\mathcal{E}_k} \ee^{V(e_k)} \overline{J}_{k}\Big]\Big),
\end{eqnarray*}
where $J_{i}:=\sum_{j=e_i}^{e_{i+1}} \ee^{V(j)-V(e_i)}$ for $i \ge
0$ and $\overline{J}_{k}:=\sum_{j=e_k}^{T^{\uparrow}(h)-1}
\ee^{V(j)-V(e_k)}$ which is well defined on $\mathcal{E}_k.$ Observe
that $\mathcal{E}_k=\{N(h)=k\}$ and recall that $N(h)$ is a
geometrical random variable with parameter $q=q(h)=P\{H \ge h\}.$
Then, the Markov property applied at times $(e_j)_{1 \le j \le k}$
yields that $\beta_2(h)$ is less or equal than
\begin{eqnarray*}
\sum_{k \ge 0}(1-q)^k q \Big(E[J_0|H_0<h] \sum_{j=0}^{k-1}
 E[\ee^{V(e_1)}|H_0<h]^j + E[\overline{J}_0|H_0\ge h]
E[\ee^{V(e_1)}|H_0<h]^k\Big),
\end{eqnarray*}
which implies that  $\beta_2(h)$ is bounded from above by
\begin{eqnarray*}
 \frac{1}{1-E[\ee^{V(e_1)}|H_0<h]} E[J_0|H_0<h]
 + \frac{q}{1-(1-q)E[\ee^{V(e_1)}|H_0<h]} E[\overline{J}_0|H_0\ge
h].
\end{eqnarray*}
Now, since $V$ is transient to $-\infty,$ then $H_0$ is almost
surely finite and $E[\ee^{V(e_1)}|H_0<h] \to E[\ee^{V(e_1)}]<1,$
when $h \to \infty.$ Recalling that $q=q(h) \to 0,$ $h \to \infty ,$
it follows that
\begin{eqnarray}
\label{eqbeta}
\beta_2(h) &\le& C\big(E[J_0|H_0<h]+q
E[\overline{J}_0|H_0\ge h]\big),
\end{eqnarray}
 for $h$ large enough.

Let us first bound $E[\overline{J}_0|H_0\ge h].$ Recall that if
$\mu$ denotes the law of $\rho_0,$ thanks to assumption $(a)$
of Theorem \ref{t:main} we can define the law
$\tilde\mu=\rho_0^\kappa \mu,$ and the law $\tilde
P=\tilde\mu^{\otimes \z}$ which is the law of a sequence of i.i.d.
random variables with law $\tilde\mu$. The definition of $\kappa$
implies that $ \int \log \rho \, \tilde \mu(d\rho) >0.$ Then, using
the explicit form of the Radon-Nykodym derivative between $P$ and $\tilde{P},$ we can write
\begin{eqnarray*}
E[\overline{J}_0|H_0\ge h] &\le& C \ee^{\kappa h}
\tilde{E}[\ee^{-\kappa V(T_h)}\overline{J}_0 {\bf 1}_{\{H_0\ge h\}}]
\\
&\le& C \tilde{E}\Big[\ee^{-\kappa(V(T_h)-h)}
\sum_{k=0}^{T_h-1}\ee^{V(k)}{\bf 1}_{\{H_0\ge h\}}\Big]
\\
&\le& C \tilde{E}\Big[\sum_{k=0}^{T_h-1}\ee^{V(k)}{\bf
1}_{\{\min_{0<k<T_h} V(k)>0\}}\Big]
\\
&\le& C \tilde{E}\Big[\sum_{k\ge0} \sum_{p=0}^{\lfloor h \rfloor}
\ee^{V(k)}{\bf 1}_{\{p\le V(k)< p+1\}}\Big]
 \\
 &\le& C  \sum_{p=0}^{\lfloor h \rfloor} \ee^{p+1}
\tilde{E}\Big[\sum_{k\ge0} {\bf 1}_{\{p\le V(k)< p+1\}}\Big].
\end{eqnarray*}
Moreover, by Markov property, we have $\tilde{E}[\sum_{k\ge0} {\bf
1}_{\{p\le V(k)< p+1\}}] \le \tilde{E}[\sum_{k\ge0} {\bf 1}_{\{0\le
V(k)< 1\}}],$ which is finite since $(V(k))_{k \ge 0}$ has a
positive drift under $\tilde{P}.$

Therefore, recalling (\ref{eqbeta}) and (\ref{iglehartthm}), we get
\begin{eqnarray}
\label{eqbeta1} \beta_2(h) \le C(E[J_0|H_0<h]+ \ee^{(1-\kappa)h})
\end{eqnarray}
and only have to bound $E[J_0|H_0<h].$ Recall that
$R=\sum_{k\geq0}\ee^{V(k)}$ and observe that $J_0 \le R.$ Moreover,
let us denote by $E^{\mathcal{I}}[\cdot]$ the expectation under
$P^{\mathcal{I}}\{\cdot\}:=P\{\cdot|\mathcal{I}\},$ with
$\mathcal{I}:=\{H=S\}.$ Then, we first observe that
$E^{\mathcal{I}}[R|H<h] \ge E[R \, {\bf 1}_{\{H=S<h\}}] \ge E[J_0
{\bf 1}_{\{H=S<h\}}].$ Furthermore, since $J_0$ depends only on
$(V(k) \, ; \, 0 \le k \le e_1)$  and since $P\{V(k) \le 0 \, ;\, k
\ge 0\}>0,$ we get, by applying the strong Markov property at time
$e_1,$ that $E[J_0 {\bf 1}_{\{H<h\}}] \le C E^{\mathcal{I}}[R|H<h],$
which implies
$$E[J_0 |H<h] \le C E^{\mathcal{I}}[R|H<h].$$
Therefore, we only have to prove that $E^{\mathcal{I}}[R|H<h] \le C
\ee^{(1-\kappa)h}.$ To this aim, we recall first that Corollary
$4.1$ in \cite{enriquez-sabot-zindy+} implies that,
$P^{\mathcal{I}}$-almost surely,
\begin{eqnarray}
\label{expo:eq5} E^{\mathcal{I}}[R|\lfloor H \rfloor]\le C
\ee^{\lfloor H \rfloor}.
\end{eqnarray}
Now, observe that $E^{\mathcal{I}}[R|H<h] \le C E^{\mathcal{I}}[R \,
{\bf 1}_{\{H<h\}}]$ and let us write
\begin{eqnarray}
E^{\mathcal{I}}[R \, {\bf 1}_{\{H<h\}}] &\le& \sum_{k=0}^{\lfloor h
\rfloor} E^{\mathcal{I}}\Big[ {\bf 1}_{\{\lfloor H \rfloor=k\}}
E^{\mathcal{I}}[R |\lfloor H \rfloor=k ]\Big] \nonumber
\\
&\le& C \sum_{k=0}^{\lfloor h \rfloor} E^{\mathcal{I}}\Big[ {\bf
1}_{\{\lfloor H \rfloor=k\}} \ee^{\lfloor H \rfloor}\Big] \nonumber
\\
 &\le& C \sum_{k=0}^{\lfloor h \rfloor}  \ee^{k} P^{\mathcal{I}}\{ \lfloor H
\rfloor=k\} \nonumber
\\
&\le& C \sum_{k=0}^{\lfloor h \rfloor}  \ee^{(1-\kappa)k}
 \le C \ee^{(1-\kappa)h},
\label{expo:eq6}
\end{eqnarray}
the second inequality is a consequence of (\ref{expo:eq5}) and the
fourth inequality due to the fact that $P^{\mathcal{I}}\{ \lfloor H
\rfloor=k\} \le c \ee^{-\kappa k}$ for some positive constant $c.$
Now assembling (\ref{expo:eq1}), (\ref{expo:eq4}), (\ref{eqbeta1})
and (\ref{expo:eq6}) concludes the proof of Lemma \ref{lemmaenv}.
 \end{proof}

\subsubsection{Important preliminary results}
\label{Dtanalogie}

Before establishing the announced lemmas, we introduce, for any $x,y
\in \z,$
$$
\tau(x,y):=\inf\{k \ge 0: \, X_{\tau(x)+k}=y\}.
$$

\noindent Recall that $A(n)=A_1(n) \cap A_2(n) \cap A_3(n) \cap A_4(n),$ where the events $(A_i(n))_{1\le i \le4}$ are defined at the beginning of Subsection \ref{themule}. Then, we have the following results.

\medskip
\begin{lemma}
 \label{l:ADT}
 Defining $DT(n) := A(n) \cap \bigcap_{j=1}^{K_n}
    \left\{ \tau(d_j,b_{j+1})< \tau(d_j,\overline{d}_j) \right\},$ we
    have
 \begin{eqnarray*}
    P\left\{ DT(n) \right\}  \to 1, \qquad n \to \infty.
     \label{ADTdef}
 \end{eqnarray*}
\end{lemma}

\begin{proof} Recalling (\ref{P:mule1}), we only have to
prove that
\begin{eqnarray}
\label{ADTeq1} E \bigg[{\bf 1}_{A(n)} \, \sum_{j=1}^{K_n}
    P_{\omega}^{d_j}\{ \tau( b_{j+1})> \tau(\overline{d}_j)\} \bigg] \to 0.
\end{eqnarray}

\noindent By (Zeitouni \cite{zeitouni}, formula (2.1.4)), we get,
for $1 \le j \le K_n$ and for all $\omega$ in $A(n):$
$$
P_{\omega}^{d_j} \left\{ \tau( b_{j+1})> \tau(\overline{d}_j) \right\} =
{\sum_{k=d_j}^{b_{j+1}-1} \ee^{V(k)} \over
\sum_{k=\overline{d}_j}^{b_{j+1}-1}\ee^{V(k)}} \le (b_{j+1}-d_j)
\ee^{V(d_j)-V(\overline{d}_j)+h_n}.
$$

\noindent Now, let us explain why $b_{K_n+1}-d_{K_n}\le 2 n$ with probability tending to
$1.$ Observe first that $b_{K_n+1}-d_{K_n}\le n + T^{\uparrow}(h_n) \circ \theta_n$ if $d_{K_n} \le n$ and  $b_{K_n+1}-d_{K_n}\le T^{\uparrow}(h_n) \circ \theta_n$ if $d_{K_n} > n.$ Therefore it is sufficient to prove that $P\{T^{\uparrow}(h_n) \ge n\} \to 0.$ But using Markov's inequality together with (\ref{expo:eq4}), we get   $P\{T^{\uparrow}(h_n) \ge n\} \le C n^{-1} \ee^{\kappa h_n} \to 0,$ when $n \to\infty.$

\noindent Moreover we have $b_{j+1}-d_j \le e_n \le C' \, n$ on $A_1(n)$ for $1
\le j \le K_n-1$ and by definition $V(d_j)-V(\overline{d}_j) \le -D_n$
for $1 \le j \le K_n.$ Therefore, we get
\begin{eqnarray*}
  E \bigg[{\bf 1}_{A(n)} \, \sum_{j=1}^{K_n}
  P_{\omega}^{d_j} \{ \tau( b_{j+1})> \tau(\overline{d}_j)\} \bigg]
  \le C \, n  E[K_n] \ee^{-D_n+h_n}.
\end{eqnarray*}

\noindent Recalling that $D_n=(1+{1 \over \kappa})\log n,$
$h_n={1-\varepsilon \over \kappa} \log n$ and since $E[K_n] \le C \,
n^{\varepsilon}$ ($K_n$ has a binomial distribution with parameter $(n,q_n)$), we obtain
\begin{eqnarray*} E \bigg[{\bf 1}_{A(n)} \,
\sum_{j=1}^{K_n}
    P_{\omega}^{d_j}\{ \tau( b_{j+1})> \tau(\overline{d}_j) \} \bigg]
     &\le& C \,  \ee^{\varepsilon(1-1/\kappa) \log n},
\end{eqnarray*}

\noindent which implies (\ref{ADTeq1}). \end{proof}

\medskip

\begin{lemma}
 \label{l:ADTgamma}
 Defining $DT^*(n):=\bigcap_{j=1}^{K_n^*}
    \left\{ \tau(b_j^*,d_j^*)< \tau(b_j^*,\gamma^*_{j})
    \right\},$ we have
 \begin{eqnarray*}
    P\{ DT^*(n)\} \to 1, \qquad n \to \infty.
     \label{ADTgammadef}
 \end{eqnarray*}
\end{lemma}

\begin{proof} Recall that $A^*(n)=\{  K_n=K_n^* \, ; \,
(a_j,b_j,c_j,d_j)= (a^*_j,b^*_j,c^*_j,d^*_j), \, 1 \le j \le K_n
\}.$ Then, let us consider
$A^{\dag}(n):=A^*(n) \cap A_3(n) \cap A_4^*(n)$ to control the
$*$-valleys, where $A_4^*(n)$ is defined by $A_4^*(n) :=
\cap_{j=1}^{K_n^*} \left\{ \gamma_{j+1}^*- a_j^* \le C'' \log n
\right\} \cap \left\{ \gamma_{1}^* \le C'' \log n \right\}.$ Using
the same arguments as in the proof of Lemma \ref{l:Aboxlength}, we
can prove that $P\{A_4^*(n)\} \to 1,$ $n \to \infty,$ for $C''$
large enough. Then, recalling that Lemma \ref{l:Aboxgap} and Lemma
\ref{model1+2} imply $P\{A^*(n) \cap A_3(n)\} \to 1,$ $n \to
\infty,$ it remains only to prove that
\begin{eqnarray}
\label{ADTgammaeq1} E \bigg[{\bf 1}_{A^{\dag}(n)} \,
\sum_{j=1}^{K_n} P_{\omega}^{b_j}\{ \tau( d_j)> \tau(\gamma^*_{j})
\} \bigg] \to 0.
\end{eqnarray}

\noindent Observe that by (Zeitouni \cite{zeitouni}, formula
(2.1.4)) we get, for $1 \le j \le K_n,$
\begin{eqnarray*}
P_{\omega}^{b_j}\{ \tau(d_j)> \tau(\gamma^*_{j}) \} &\le& (d_j-b_j)
\ee^{H^{(j)}-(V(\gamma^*_{j})-V(b_j))}
\\
&\le& C \log n \, \ee^{H^{(j)}-(V(\gamma^*_{j})-V(b_j))},
\end{eqnarray*}

\noindent the second inequality being a consequence of $\omega \in
A^*(n)\cap A_4^*(n).$ Then, to bound
$\ee^{H^{(j)}-(V(\gamma^*_{j})-V(b_j))}$ from above, observe that
(\ref{iglehartthm}) implies $ P\{ H^{(j)}> {(1+\varepsilon') \over
\kappa} \log n\} \sim
n^{-(\varepsilon'+\varepsilon)}=o(n^{-\varepsilon}),$ $n \to
\infty,$ for any $\varepsilon'>0,$ which yields that
$P\{\bigcap_{j=1}^{K_n} \{H^{(j)} < {(1+\varepsilon') \over \kappa}
\log n \}\}$ tends to $1$, when $n$ tends to $\infty.$ Therefore,
recalling (\ref{ADTgammaeq1}), we only have to prove that
\begin{eqnarray}
\label{ADTgammaeq2} C (\log n) n^{{(1+\varepsilon') \over \kappa}} E
\bigg[{\bf 1}_{A^{\dag}(n)} \, \sum_{j=1}^{K_n}
\ee^{-(V(\gamma^*_{j})-V(b_j))} \bigg] \to 0.
\end{eqnarray}

\noindent Since $\gamma_{j}^*- b_{j-1} \le C'' \log n$ on $A_4^*(n)$
and $b_{j}-b_{j-1} \ge n^{1-3\varepsilon}$ on $A_3(n),$ we get
$b_j-\gamma_{j}^*\ge{1 \over 2}n^{1-3\varepsilon}$ for $2 \le j \le
K_n$ on $A^{\dag}(n)$, for all large $n.$ Similarly,
$\gamma_{0}^*\le C'' \log n$ on $A_4^*(n)$ and $b_1\ge
n^{1-3\varepsilon}$ on $A_3(n)$ yield $b_1-\gamma_{1}^*\ge{1 \over
2}n^{1-3\varepsilon}$ on $A^{\dag}(n).$
 Therefore, recalling the definition of $b_j,$ we can use Lemma \ref{lemmeLD} and obtain\begin{eqnarray*}
 P\{ A^{\dag}(n) \, ;\, V(b_j)- V(\gamma^*_{j}) \ge
-n^{1-3\varepsilon \over 2}\} &\le& P\{ V({1 \over
2}n^{1-3\varepsilon}) \ge -n^{1-3\varepsilon \over 2}\}
\\
&\le& \ee^{-{n^{1-3\varepsilon} \over 2} I\big(-2 n^{-{1-3\varepsilon
\over 2}}\big)}  =o(n^{-\varepsilon}),
\end{eqnarray*}

\noindent for any $1 \le j \le K_n,$ since $I(0)>0.$ This result implies that the
term on the left-hand side in (\ref{ADTgammaeq2}) is bounded from
above by $C \log n \ n^{{(1+\varepsilon') \over \kappa}} E[K_n]
\ee^{-{n^{1-3\varepsilon} \over 2}}.$ Then, since $E[K_n] \le C \,
n^{\varepsilon},$ this concludes the proof of Lemma
\ref{l:ADTgamma}. \end{proof}

\medskip

\begin{lemma}
 \label{l:IA}
 For any $0<\eta<\varepsilon({1 \over \kappa}-1),$ let us introduce the following event
  $IA(n) :=A(n) \cap \left\{  \sum_{j=1}^{K_n} \tau(d_j,b_{j+1})<n^{1/\kappa-\eta}  \right\}.$
   Then, we have
 \begin{eqnarray*}
    P\{ IA(n)\}  \to 1, \qquad n \to \infty.
     \label{AIAdef}
 \end{eqnarray*}
\end{lemma}

\begin{proof} Recalling that $P\{ K_n \ge 2 C_I
n^{\varepsilon}\} \to 0,$ $n \to \infty,$ and that Lemma \ref{l:ADT}
implies that $P\{ DT(n)\} \to 1,$ $n \to \infty,$ it only remains to
prove
\begin{eqnarray*}
 \p \bigg\{DT(n) \cap \bigg\{ \sum_{j=1}^{\lfloor 2 C_I n^{\varepsilon} \rfloor}
  \tau(d_j,b_{j+1})\ge n^{1/\kappa-\eta} \bigg\}  \bigg\} \to 0, \qquad n \to \infty.
 \end{eqnarray*}

\noindent Using Markov inequality, we have to prove that
\begin{eqnarray}
 \e \bigg[ {\bf 1}_{DT(n)} \, \sum_{j=1}^{\lfloor 2 C_I n^{\varepsilon} \rfloor}
  \tau(d_j,b_{j+1}) \bigg]= o \left({1 \over n^{1/\kappa-\eta}}\right),
   \qquad n \to \infty.
   \label{AIAeq1}
 \end{eqnarray}

\noindent Furthermore, by definition of the event $DT$ (see Lemma
\ref{l:ADT}), we get
\begin{eqnarray*}
 \e\bigg[ {\bf 1}_{DT(n)} \, \sum_{j=1}^{\lfloor 2 C_I n^{\varepsilon} \rfloor}
  \tau(d_j,b_{j+1}) \bigg] &\le& E \bigg[{\bf 1}_{A(n)} \, \sum_{j=1}^{\lfloor 2 C_I n^{\varepsilon} \rfloor}
   E_{\omega, |\overline{d}_j}^{d_j} [\tau(b_{j+1})] \bigg]
  \\
  &\le& E \bigg[ {\bf 1}_{A(n)} \, \sum_{j=1}^{\lfloor 2 C_I n^{\varepsilon} \rfloor}
   E_{\omega, |\overline{d}_j}^{\overline{d}_j}[\tau(b_{j+1}) ] \bigg].
    \end{eqnarray*}

\noindent Applying successively the strong Markov property at
$\overline{d}_{\lfloor 2 C_I n^{\varepsilon} \rfloor},$ \dots,
$\overline{d}_{2},$ $\overline{d}_1,$ this implies
\begin{eqnarray*}
 \e\bigg[ {\bf 1}_{DT(n)} \, \sum_{j=1}^{\lfloor 2 C_I n^{\varepsilon} \rfloor}
  \tau(d_j,b_{j+1})\bigg] &\le& 2 C_I n^{\varepsilon} \e_{|0}
 [\tau(T^{\uparrow}(h_n)-1)].
\end{eqnarray*}

\noindent Therefore, Lemma \ref{lemmaenv} implies
\begin{eqnarray*}
 \e \bigg[ {\bf 1}_{DT(n)} \, \sum_{j=1}^{\lfloor 2 C_I n^{\varepsilon} \rfloor}
  \tau(d_j,b_{j+1}) \bigg] \le  C  n^{\varepsilon}  \ee^{h_n}
\le C n^{{1 \over \kappa}-\varepsilon \left({1 \over
  \kappa}-1 \right)},
\end{eqnarray*}

\noindent which yields (\ref{AIAeq1}) and concludes the proof, since
$0<\eta<\varepsilon({1 \over \kappa}-1).$ \end{proof}

\medskip

\begin{lemma}
 \label{l:valley-*valley} We have
  \begin{eqnarray*}
P\{(a_1,b_1,c_1,d_1) \neq
(a_1^*,b_1^*,c_1^*,d_1^*) \} = o(n^{-\varepsilon}), \qquad n \to \infty.
 \end{eqnarray*}
\end{lemma}

\begin{proof} Since $\gamma_1^*$ is a negative record for the potential $V,$ it is sufficient to prove that there is no excursion higher than $h_n$ before $\gamma_1^*.$
In a first step, we prove that for $C$ large enough
  \begin{eqnarray}
  \label{eq:valley*valley1}
P\{  \gamma_1^*  \ge C  \log n  \} = o(n^{-\varepsilon}), \qquad n \to \infty.
 \end{eqnarray}
Indeed, applying Lemma \ref{lemmeLD}, we get
  \begin{eqnarray*}
P\{  \gamma_1^*  \ge C  \log n  \} &\le& P\{ V(C \log n)  \ge -D_n \}
\\
 &\le& \exp \big\{-C I({1+\kappa^{-1} \over C}) \log n\big\}  = o(n^{-\varepsilon}),
 \end{eqnarray*}
by choosing $C$ so that $C I({1+\kappa^{-1} \over C})> \varepsilon,$ which is possible since $I(0)>0.$ 

In a second step, we prove that the probability that there is an excursion higher than $h_n$ before $C \log n$ is a $o(n^{-\varepsilon}).$ Since the number of excursions before $C \log n$ is bounded by $C \log n,$ we will prove that
  \begin{eqnarray}
  \label{eq:valley*valley2}
P\left\{ \max_{0\le i \le C \log n} H_i  \ge h_n \right\} = o(n^{-\varepsilon}), \qquad n \to \infty.
 \end{eqnarray}
 But this result is obvious. Indeed, using (\ref{iglehartthm}) we obtain that the probability term in (\ref{eq:valley*valley2}) is less than $C \log n \, \ee^{-\kappa h_n}= o(n^{-\varepsilon}).$
Now assembling (\ref{eq:valley*valley1}) and (\ref{eq:valley*valley2}) concludes the proof of Lemma \ref{l:valley-*valley}.

 \end{proof}

\subsection{Proof of Proposition \ref{P:prop1}}
Since the time spent on $\z_-$ is almost surely finite, we reduce
our study to the random walk in random environment reflected at $0$
and observe that
\begin{eqnarray*}
     \e \left[ \ee^{- \lambda_n \, \tau(e_n)} \right]= \e_{|0}\left[ \ee^{- \lambda_n \, \tau(e_n)}
     \right]+o(1), \qquad n \to \infty,
\end{eqnarray*}

\noindent where $\e_{|0}$ denotes the expectation under the law
$\p_{|0}$ of the random walk in the random environment $\omega$
(under $P$) reflected at $0.$

Furthermore, by definition, $\tau(e_n)$ satisfies
$$
\tau(b_1)+\sum_{j=1}^{K_n-1}\{\tau(b_j,d_{j}) + \tau(d_j,b_{j+1})\}
\le \tau(e_n) \le \tau(b_1)+\sum_{j=1}^{K_n} \{\tau(b_j,d_{j}) +
\tau(d_j,b_{j+1})\},
$$

\noindent such that we easily get that $\e_{|0} \left[ \ee^{-
\lambda_n \, \tau(e_n)} \right]$ belongs to
\begin{eqnarray*}
     \left[\e_{|0} \left[ \ee^{- \lambda_n \, \left( \tau(b_1)+\sum_{j=1}^{K_n}
     \{\tau(b_j,d_{j}) + \tau(d_j,b_{j+1})\}\right)} \right] , \,
      \e_{|0} \left[ \ee^{- \lambda_n \, \left( \tau(b_1)+\sum_{j=1}^{K_n-1}
      \{\tau(b_j,d_{j}) + \tau(d_j,b_{j+1})\} \right)} \right]  \right].
\end{eqnarray*}

\noindent  Let us first recall that Lemma \ref{l:ADT} and Lemma
\ref{l:IA} imply that $P\{DT(n) \cap IA(n)\} \to 1,$ $n \to \infty.$
Then, we get that the lower bound in the previous interval is equal
to
\begin{eqnarray*}
    && \e_{|0} \left[{\bf 1}_{DT(n) \cap IA(n)} \ee^{- \lambda_n
    \left( \tau(b_1)+\sum_{j=1}^{K_n}
     \{\tau(b_j,d_{j}) + \tau(d_j,b_{j+1})\} \right)} \right] +o(1)
    \\
    &=&
    \e_{|0} \left[ {\bf 1}_{DT(n) \cap IA(n)} \, \ee^{- \lambda_n
    \sum_{j=1}^{K_n} \tau(b_j, d_j)} \right]+o(1)
    \\
   &=& \e_{|0} \left[  \ee^{- \lambda_n
    \sum_{j=1}^{K_n} \tau(b_j, d_j)}  \right]+o(1).
\end{eqnarray*}

\noindent Then, applying the strong Markov property for the random
walk successively at $\tau(b_{K_n}),$ $\tau(b_{K_n-1}),$ \dots,
$\tau(b_{2})$ and $\tau(b_{1})$ we get
\begin{eqnarray*}
  \e_{|0} \left[  \ee^{- \lambda_n
    \sum_{j=1}^{K_n} \tau(b_j, d_j)}  \right] &=&
    E \bigg[\prod_{j=1}^{K_n}  E_{\omega,|0}^{b_j} \left[ \ee^{- \lambda_n
     \tau(d_j)} \right] \bigg]
    \\
  &=& E \bigg[{\bf 1}_{A^*(n)} \, \prod_{j=1}^{K_n^*}  E_{\omega,|0}^{b_j^*} \left[ \ee^{- \lambda_n
     \tau(d_j^*)} \right] \bigg]+o(1)
     \\
   &=& E \bigg[ \prod_{j=1}^{K_n^*}  E_{\omega,|0}^{b_j^*} \left[ \ee^{- \lambda_n
     \tau(d_j^*)} \right] \bigg]+o(1),
\end{eqnarray*}

\noindent the second equality being a consequence of Lemma
\ref{model1+2}. Then, since Lemma \ref{l:ADTgamma} implies
$\p\{DT^*(n) \} \to 1,$ we have
\begin{eqnarray*}
    \e_{|0} \left[  \ee^{- \lambda_n
    \sum_{j=1}^{K_n} \tau(b_j, d_j)}  \right] &=&
    E \bigg[\prod_{j=1}^{K_n^*}  E_{\omega,|0}^{b_j^*} \left[{\bf 1}_{DT^*(n)} \,
     \ee^{- \lambda_n \tau(d_j^*)} \right] \bigg]+o(1)
    \\
   &=& E \bigg[\prod_{j=1}^{K_n^*}  E_{\omega,|\gamma_{j}^*}^{b_j^*} \left[{\bf 1}_{DT^*(n)}
     \, \ee^{- \lambda_n \tau(d_j^*)} \right]\bigg]+o(1)
     \\
   &=& E \bigg[\prod_{j=1}^{K_n^*}  E_{\omega,|\gamma_{j}^*}^{b_j^*} \left[
     \ee^{- \lambda_n \tau(d_j^*)} \right] \bigg]+o(1),
\end{eqnarray*}

\noindent Since $\p\{K_n=K_n^* \} \to 1,$ and $\p\{ K_n \le
\overline{K}_n \} \to 1,$ with  $\overline{K}_n= \lceil n q_n
(1+n^{-\varepsilon/4}) \rceil,$ we get
\begin{eqnarray*}
  \e_{|0} \left[ \ee^{- \lambda_n \, \tau(e_n)} \right] \ge
     E \bigg[ \prod_{j=1}^{\overline{K}_n}  E_{\omega,|\gamma_{j}^*}^{b_j^*} \left[
     \ee^{- \lambda_n \tau(d_j^*)} \right]\bigg]
     +o(1).
\end{eqnarray*}

\noindent Then, applying the strong Markov property (for the
potential $V$) successively at times $\gamma_{\overline{K}_n}^*,$
\dots, $\gamma_{2}^*$ and observing that the
$\left(E_{\omega,|\gamma_{j}^*}^{b_j^*} \left[ \ee^{- \lambda_n
\tau(d_j^*)} \right]\right)_{1 \le j \le \overline{K}_n}$ are i.i.d.
random variables, we obtain that
\begin{eqnarray*}
  \e_{|0} \left[ \ee^{- \lambda_n \, \tau(e_n)} \right] \ge
     E \left[ E_{\omega,|\gamma_{1}^*}^{b_1^*}
     \ee^{- \lambda_n \tau(d_1^*)}
     \right]^{\overline{K}_n}+o(1).
\end{eqnarray*}

\noindent Using Lemma \ref{l:valley-*valley} and recalling that
$\overline{K}_n= \lceil n q_n
(1+n^{-\varepsilon/4}) \rceil=O(n^{\varepsilon}),$ $n \to \infty,$ the strong
Markov property applied at $\gamma_{1}^*$ yields
\begin{eqnarray*}
  \e_{|0} \left[ \ee^{- \lambda_n \, \tau(e_n)} \right] \ge
     E \left[ E_{\omega,|0}^{b_1} \left[
     \ee^{- \lambda_n \tau(d_1)} \right]
     \right]^{\overline{K}_n}+o(1).
\end{eqnarray*}

\noindent Using similar arguments for the upper bound in the
aforementioned interval, we get
\begin{eqnarray*}
  \e_{|0} \left[ \ee^{- \lambda_n \, \tau(e_n)} \right] \in \left[
     E \left[ E_{\omega,|0}^{b_1} \left[
     \ee^{- \lambda_n \tau(d_1)} \right]
     \right]^{\overline{K}_n}+o(1) \, , \,
     E \left[ E_{\omega,|0}^{b_1} \left[
     \ee^{- \lambda_n \tau(d_1)} \right]
     \right]^{\underline{K}_n} +o(1) \right].
\end{eqnarray*}
with $\underline{K}_n:= \lfloor n q_n
(1-n^{-\varepsilon/4}) \rfloor.$ Furthermore, observe that we have $E \left[ E_{\omega,|0}^{b_1}
\left[ \ee^{- \lambda_n \tau(d_1)} \right] \right]=E \left[
E_{\omega,|a_1}^{b_1} \left[ \ee^{- \lambda_n \tau(d_1)} \right]
\right]+o(n^{-\varepsilon}).$ This is a consequence of Lemma
\ref{l:Aboxlength}, definition of $a$ and the fact that
(\ref{iglehartthm}) implies $ P\{ H^{(1)}
> {(1+\varepsilon') \over \kappa} \log n\} \sim
 n^{-(\varepsilon'+\varepsilon)}=o(n^{-\varepsilon}),$ $n \to
\infty,$ for any $\varepsilon'>0,$ which gives
$$
E \left[ P_{\omega}^{b_1} \left\{ \tau(a_1) < \tau(d_1) \right\}
\right] \le C \log n \, \ee^{ {(1+\varepsilon') \over \kappa} \log n
- D_n} = o(n^{-\varepsilon}).
$$

\noindent This concludes the proof of Proposition \ref{P:prop1}.
\hfill$\Box$

\section{Annealed Laplace transform for the exit time from a deep valley}
\label{sec:part2}

This section is devoted to the proof of the linearization. It
involves $h$-processes theory and ``sculpture'' of a typical deep
valley. To ease notations, we shall use $a,$ $b$, $c,$ and $d$
instead of $a_1,$ $b_1,$ $c_1$ and $d_1.$ Moreover, let us
introduce, for any random variable $Z\ge0,$ the functional
\begin{eqnarray}
\label{RnlambdaZ} R_n(\lambda,Z):= E\bigg[{1 \over
1+\frac{\lambda}{n^{1/\kappa}} Z}\bigg],
\end{eqnarray}
and the two important random variables given by
\begin{eqnarray}
\widehat{M}_1&:=&\sum_{x=a+1}^{d-1}
\ee^{-(\widehat{V}(x)-\widehat{V}(b))},
\\
M_2&:=&\sum_{x=b}^{d-1}\ee^{V(x)-V(c)},
\end{eqnarray}
where $\widehat{V}$ is defined below in (\ref{def:Vchapeau}).
Then, the result can be expressed in the following way.
\medskip

\begin{proposition}
 \label{P:prop2} For any $\xi>0,$ we have, for all large $n,$
 \begin{eqnarray*}
 R_n(\ee^{\xi} \lambda,2  \ee^{H^{(1)}}\widehat{M}_1 M_2)\! +\!o(n^{-\varepsilon})\!
\le \! E\big[E_{\omega,|a}^{b} [\ee^{- \lambda_n \tau(d)}] \big]\!
\le \! R_n(\ee^{-\xi} \lambda,2  \ee^{H^{(1)}}\widehat{M}_1 M_2)
\!+\!o(n^{-\varepsilon}).
\end{eqnarray*}
\end{proposition}

\medskip
\subsection{Two $h$-processes} \label{h-potential}
 In order to estimate $E_{\omega,|a}^{b} \left[ \ee^{-
\lambda_n \tau(d)} \right],$ we decompose the passage from $b$ to
$d$ into the sum of a random geometrically distributed number,
denoted by $N$, of unsuccessful attempts to reach $d$ from $b$ (i.e. excursions of the particle  from $b$ to $b$ which do not hit $d$), followed by a successful attempt.   More precisely, since $N$ is a
geometrically distributed random variable with parameter $1-p$
satisfying (see \cite{zeitouni}, formula (2.1.4))
\begin{eqnarray}
\label{1-p} 1-p &=&  \omega_b \, {\ee^{V(b)} \over \sum_{x=b}^{d-1}
\ee^{V(x)}},
\end{eqnarray}

\noindent we can write $ \tau(d)=\sum_{i=1}^{N} F_i +G,$ where the
$F_i$'s are the successive i.i.d. failures and $G$ the first
success. The accurate estimation of the time spent by each
(successful and unsuccessful) attempt leads us to consider two
$h$-processes where the random walker evolves in two modified
potentials, one corresponding to the conditioning on a failure (see
the potential $\widehat{V}$ and Lemma \ref{lemmafailure}) and the
other to the conditioning on a success (see the potential $\bar{V}$
and Lemma \ref{lemmasuccess}).

\subsubsection{The failure case: the $h$-potential $\widehat{V}$}
\label{Vcheck}

Let us fix a realization of $\omega.$ To introduce the $h$-potential
$\widehat{V},$ we consider the valley $a<b<c<d$ and define $h(x):=
P_{\omega}^{x} \{\tau(b) < \tau(d)\}.$ For any $b<x<d,$
we introduce $\widehat{\omega}_x:=\omega_x {h(x+1)\over h(x)}.$ Since $h$ is a harmonic function, we have $1-\widehat{\omega}_x=(1-\omega_x) {h(x-1)\over h(x)}.$ Now, $\widehat{V}$ can be defined for $x \ge b$  by
\begin{eqnarray}
\label{def:Vchapeau}
\widehat{V}(x) := V(b)+ \sum_{i=b+1}^x \log {1-\widehat{\omega}_i \over\widehat{\omega}_i}.
\end{eqnarray}
We obtain for any $b\le x<y<d,$
\begin{eqnarray}
\label{Vcheckeq0}
\widehat{V}(y)-\widehat{V}(x)=\left(V(y)-V(x)\right)+ \log \bigg({
h(x) \, h(x+1) \over h(y) \, h(y+1)} \bigg).
\end{eqnarray}

\noindent Since $h(x)$ is a decreasing function of $x$ by definition, we get
\begin{eqnarray}
\label{Vcheckeq1}
 { h(x) \, h(x+1) \over h(y) \, h(y+1) } \ge 1.
\end{eqnarray}

\noindent Thus we obtain for any $b\le x < y \le c,$
\begin{equation}\label{Vcheckeq2}
\widehat{V}(y)-\widehat{V}(x) \ge V(y)-V(x).
\end{equation}

\medskip

\begin{lemma} For any environment $\omega,$ we have
 \label{lemmafailure}
\begin{equation}\label{expF}
   E_{\omega} \left[ F_1 \right]= 2 \, \omega_b \bigg( \sum_{i=a+1}^{b-1}
    \ee^{-(V(i)-V(b))} + \sum_{i=b}^{d-1}
    \ee^{-(\widehat{V}(i)-\widehat{V}(b))}\bigg),
\end{equation}
\noindent and
\begin{equation}\label{expF^2}
    E_{\omega} \left[ F_1^2 \right]= 4 \omega_b \, R^+ + 4(1-
    \omega_b) \, R^-,
\end{equation}
\noindent where
\begin{eqnarray*}
  R^+ &:=& \sum_{i=b+1}^{d-1} \bigg(1+2 \sum_{j=b}^{i-2} \ee^{\widehat{V}(j)-\widehat{V}(i-1
  )} \bigg) \bigg( \ee^{-(\widehat{V}(i-1)-\widehat{V}(b))} +2 \sum_{j=i+1}^{d-1}
  \ee^{-(\widehat{V}(j-1)-\widehat{V}(b))} \bigg),
   \\
  R^- &:=& \sum_{i=a+1}^{b-1} \bigg(1+2 \sum_{j=i+2}^{b} \ee^{V(j)-V(i+1
  )} \bigg)\bigg( \ee^{-(V(i+1)-V(b))} +2 \sum_{j=a+1}^{i-1}
  \ee^{-(V(j+1)-V(b))}\bigg).
\end{eqnarray*}
\end{lemma}

\medskip

\begin{remark}
\label{Goldsheid-alili} Alili \cite{alili} and Goldsheid
\cite{goldsheid} prove a similar result for a non-conditioned
hitting time. Here we give the proof in order to be self-contained.
\end{remark}

\medskip

\begin{proof} Let us first introduce
\begin{eqnarray*}
% \nonumber to remove numbering (before each equation)
  N_i^+&:=&\sharp\{k < \tau(b): X_k=i-1, X_{k+1}=i \}, \qquad i>b,\\
  N_i^-&:=&\sharp\{k < \tau(b): X_k=i+1, X_{k+1}=i \},  \qquad i<b,
\end{eqnarray*}
and the quenched probability in the environment $\widehat{\omega},$ denoted by $P_{\widehat{\omega}}.$ Then, observe that, under $P_{\widehat{\omega}},$ for $i>b$ and
conditionally on $N_{i}^+=x$, $N_{i+1}^+$ is the sum of $x$
independent geometrical random variables with parameter
$\widehat{\omega}_{i} \in (0,1).$ It means that $E_{\widehat{\omega}} [
 N_{i+1}^+ |N_{i}^+=x]={x \over \widehat{\rho}_{i}}$ and ${\rm Var}_{\widehat{\omega}} [
 N_{i+1}^+ |N_{i}^+=x]={x \over\widehat{\omega}_{i} \widehat{\rho}_{i}^2}.$
Similarly, under $P_{\omega},$ for $i<b$ and conditionally on
$N_{i}^-=x$, $N_{i-1}^-$ is the sum of $x$ independent geometrical
random variables with parameter $1-{\omega}_{i}.$ It means that
$E_{{\omega}} [
 N_{i-1}^- |N_{i}^-=x]={x{\rho}_{i}}$ and ${\rm Var}_{{\omega}} [
 N_{i-1}^- |N_{i}^-=x]={x{\rho}_{i}^2 \over(1-{\omega}_{i})}.$

\noindent  Since
$$E_{\omega} [ F_1 ]= 2 \omega_b \, E_{\widehat{\omega}}
[ \sum_{b+1}^{d-1} N_i^+ ] + 2(1- \omega_b) \, E_{\omega} [
\sum_{a+1}^{b-1} N_i^- ],$$

\noindent an easy calculation yields (\ref{expF}).

To calculate $E_{\omega}[ F_1^2],$ observe first that
$$E_{\omega}[ F_1^2 ]= 4 \omega_b \, E_{\widehat{\omega}}
\bigg[ (\sum_{i=b+1}^{d-1} N_i^+)^2 \bigg] + 4(1- \omega_b) \,
E_{\omega} \bigg[(\sum_{i=a+1}^{b-1} N_i^-)^2\bigg].$$

\noindent Then, it remains to prove that $E_{\widehat{\omega}} [
(\sum_{b+1}^{d-1} N_i^+)^2 ]=R^+$ and $E_{\omega} [
(\sum_{a+1}^{b-1} N_i^-)^2 ]=R^-.$ We will only treat
$E_{\widehat{\omega}} [ (\sum_{b+1}^{d-1} N_i^+)^2 ]$, the case of
$E_{\omega} [ (\sum_{a+1}^{b-1} N_i^-)^2 ]$ being similar. We get
first
\begin{eqnarray}
\label{defR_0}
  E_{\widehat{\omega}} \bigg[ (\sum_{b+1}^{d-1} N_i^+)^2 \bigg]=
  \sum_{i=b+1}^{d-1} E_{\widehat{\omega}}[ ( N_i^+)^2] +
2 \sum_{i=b+1}^{d-1} \sum_{j=i+1}^{d-1} E_{\widehat{\omega}}[ N_i^+
N_j^+].
\end{eqnarray}

\noindent Observe that $E_{\widehat{\omega}}\left[ N_i^+ N_j^+
\right]=E_{\widehat{\omega}}\left[ N_i^+ E_{\widehat{\omega}}\left[
N_j^+ \, | \, N_i^+,\dots ,N_{j-1}^+\right] \right]=
E_{\widehat{\omega}}\left[ N_i^+ {N_{j-1}^+ \over
\widehat{\rho}_{j-1}} \right],$ for $i<j,$ so that we get, by
iterating,
$$
E_{\widehat{\omega}}\left[ N_i^+ N_j^+
\right]=E_{\widehat{\omega}}\left[ (N_i^+)^2 \right]{1 \over
\widehat{\rho}_{j-1} \dots \widehat{\rho}_{i} }.
$$

\noindent Recalling (\ref{defR_0}), this yields
\begin{eqnarray}
 \nonumber
  E_{\widehat{\omega}} \bigg[ (\sum_{b+1}^{d-1} N_i^+)^2 \bigg]&=&\sum_{i=b+1}^{d-1}
   E_{\widehat{\omega}}\left[ ( N_i^+)^2
   \right] \bigg(1+2 \sum_{j=i+1}^{d-1} {1 \over \widehat{\rho}_{i}
   \dots \widehat{\rho}_{j-1}} \bigg)
 \\
   &=& \sum_{i=b+1}^{d-1} E_{\widehat{\omega}}\left[ ( N_i^+)^2 \right]
    \bigg(1+2 \sum_{j=i+1}^{d-1} \ee^{-(\widehat{V}(j-1)-\widehat{V}(i-1
  ))} \bigg).
  \label{R0eq1}
\end{eqnarray}

\noindent Now, observe that $E_{\widehat{\omega}}\left[ (N_i^+)^2
\right]=E_{\widehat{\omega}}\left[  E_{\widehat{\omega}}\left[
(N_i^+)^2 | N_{i-1}^+\right] \right],$ which implies
\begin{eqnarray*}
E_{\widehat{\omega}}\left[ (N_i^+)^2 \right]=
E_{\widehat{\omega}}\bigg[ \sum_{k \ge 1} E_{\widehat{\omega}}[
G_1^{(i)}+ \dots + G_k^{(i)} ] {\bf 1}_{\{N_{i-1}^+=k\}} \bigg].
\end{eqnarray*}

\noindent Since the $G_{\cdot}^{(i)}$'s are i.i.d., we get $
E_{\widehat{\omega}}[ G_1^{(i)}+ \dots + G_k^{(i)}]=k {\rm
Var}_{\widehat{\omega}}[ G_1^{(i)}]+k^2 E_{\widehat{\omega}}[
G_1^{(i)}]^2.$ Recalling that $E_{\widehat{\omega}} [ G_1^{(i)}]={1
\over \widehat{\rho}_{i-1}}$ and ${\rm Var}_{\widehat{\omega}}
[G_1^{(i)}]={1 \over\widehat{\omega}_{i-1} \widehat{\rho}_{i-1}^2},$
this yields
\begin{eqnarray}
\nonumber E_{\widehat{\omega}}\left[ (N_i^+)^2 \right]&=&
{E_{\widehat{\omega}}\left[ N_{i-1}^+ \right]
\over\widehat{\omega}_{i-1} \widehat{\rho}_{i-1}^2}+
{E_{\widehat{\omega}}\left[ (N_{i-1}^+)^2 \right] \over
\widehat{\rho}_{i-1}^2}
 \\
&=& {1 \over \widehat{\omega}_{i-1}\widehat{\rho}_{b+1} \dots
\widehat{\rho}_{i-2} \widehat{\rho}_{i-1}^2}+
{E_{\widehat{\omega}}\left[ (N_{i-1}^+)^2 \right] \over
\widehat{\rho}_{i-1}^2}. \label{Ni^2}
\end{eqnarray}

\noindent Denoting $W_{b+1}:=1$ and $W_{i}:=(\widehat{\rho}_{b+1}
\dots \widehat{\rho}_{i-1})^2 E_{\widehat{\omega}}\left[ (N_i^+)^2
\right]$ for $b+1<i<d,$ (\ref{Ni^2}) becomes
\begin{eqnarray*}
  W_{i}-W_{i-1} &=&  {\widehat{\rho}_{b+1}
\dots \widehat{\rho}_{i-1} \over \widehat{\omega}_{i-1}}=
\widehat{\rho}_{b+1} \dots \widehat{\rho}_{i-1} +
\widehat{\rho}_{b+1} \dots \widehat{\rho}_{i-2},
\end{eqnarray*}

\noindent the second equality being a consequence of
${1/\widehat{\omega}_{i-1}}=\widehat{\rho}_{i-1}+1.$ Therefore, we
have $W_{i} = \sum_{b+2}^{i}(W_{j}-W_{j-1}) + W_{b+1} =
\widehat{\rho}_{b+1} \dots \widehat{\rho}_{i-1}+ 2(1+
\sum_{b+1}^{i-2} \widehat{\rho}_{b+1} \dots \widehat{\rho}_{j}), $
which implies
\begin{eqnarray}
\nonumber E_{\widehat{\omega}}\left[ (N_i^+)^2 \right]&=& {1 \over
\widehat{\rho}_{b+1} \dots \widehat{\rho}_{i-1}} + 2
   \sum_{j=b}^{i-2} {\widehat{\rho}_{b+1} \dots \widehat{\rho}_{j} \over
   (\widehat{\rho}_{b+1} \dots \widehat{\rho}_{i-1})^2}
 \\
&=& \ee^{-(\widehat{V}(i-1)-\widehat{V}(b))} + 2
   \sum_{j=b}^{i-2} \ee^{\widehat{V}(j)-2\widehat{V}(i-1)+\widehat{V}(b)}.
   \label{Ni^2+}
\end{eqnarray}

\noindent Assembling (\ref{R0eq1}) and (\ref{Ni^2+}) yields
(\ref{expF^2}).\end{proof}

\subsubsection{The success case: the $h$-potential $\bar{V}$}
\label{Vbar}

 In a similar way, we introduce the $h$-potential $\bar{V}$
by considering the valley $a<b<c<d$ and defining $g(x):=
P_{\omega}^{x} \{\tau(d) < \tau(b)\}.$ For any $b<x<d,$
we introduce $\bar{\omega}_x:=\omega_x {g(x+1)\over g(x)}.$ Since $g$ is a harmonic function, we have $1-\bar{\omega}_x=(1-\omega_x) {g(x-1)\over g(x)}.$  Then, $\bar{V}$ can be defined for $x \ge b$  by
$$
\bar{V}(x) := V(b)+ \sum_{i=b+1}^x \log {1-\bar{\omega}_i \over\bar{\omega}_i}. $$
We have the following result for any $b<x<y\le d,$ 
\begin{eqnarray}
\label{Vbareq0} \bar{V}(y)-\bar{V}(x)=\left(V(y)-V(x)\right)+ \log
\bigg({ g(x) \, g(x+1) \over g(y) g(y+1)} \bigg).
\end{eqnarray}

\noindent Since $g(x)$ is a increasing function of $x$ by definition, we get
\begin{eqnarray}
\label{Vbareq1} { g(x) \, g(x+1) \over g(y) \, g(y+1)}\le 1.
\end{eqnarray}

\noindent Therefore, we obtain for any $c\le x < y \le d,$
\begin{eqnarray}
\label{Vbareq2}
   \bar{V}(y)-\bar{V}(x) \le V(y)-V(x).
\end{eqnarray}

Using the same arguments as in the failure case, we get the
following result.

\medskip

\begin{lemma} For any environment $\omega,$ we have
 \label{lemmasuccess}
\begin{equation}\label{expS}
   E_{\omega}[G] \le1 + \sum_{i=b+1}^{d}
    \sum_{j=i}^{d} \ee^{\bar{V}(j)-\bar{V}(i)}.
\end{equation}
\end{lemma}

\medskip

\subsection{Preparatory lemmas}

The study of a typical deep valley involves the following event
\begin{eqnarray*}
  A_5(n) &:=& \left\{ \max\{V^{\uparrow}(a,b)\, ; \, -V^{\downarrow}
  (b,c) \, ; \, V^{\uparrow}(c,d)\}
  \le \delta \log n  \right\},
\end{eqnarray*}

\noindent where $\delta > \varepsilon /\kappa.$ In words, $A_5(n)$
ensures that the potential does not have excessive fluctuations in a
typical box. Moreover, we have the following result.
\medskip

\begin{lemma}
 \label{l:Apotfluctu}
 For any $\delta > \varepsilon / \kappa,$
 \begin{eqnarray*}
   P\{ A_5(n)\} =1- o( n^{-\varepsilon}), \qquad n \to \infty.
        \label{Apotfluctudef}
 \end{eqnarray*}
\end{lemma}

\medskip

\begin{proof} We easily observe that the proof of Lemma
\ref{l:Apotfluctu}
 boils down to showing that
\begin{eqnarray}
  P\{  V^{\uparrow}(a,b) \ge \delta \log n \}
  &=& o( n^{-\varepsilon}), \qquad n \to
  \infty,
  \label{Apotfluctueq1}
  \\
  P\{-V^{\downarrow}(b,c) \ge \delta \log n\}
  &=& o( n^{-\varepsilon}), \qquad n \to
  \infty,
  \label{Apotfluctueq2}
  \\
  P\{V^{\uparrow}(c,d) \ge \delta \log n\}
  &=& o( n^{-\varepsilon}),
\qquad n \to \infty.
  \label{Apotfluctueq3}
\end{eqnarray}

In order to prove (\ref{Apotfluctueq3}), let us first observe the
following trivial inequality
$$P\{V^{\uparrow}(c,d) \ge \delta \log n
\} \le P\{V^{\uparrow}(T_{1}^{\uparrow},d) \ge \delta \log n \}.$$

\noindent Looking at the proof of (\ref{Aboxlengtheq1}), we observe
that $P\{ d-T_{1}^{\uparrow} \ge C \log n\} = o(
n^{-\varepsilon'}),$ for any $\varepsilon'>0,$ by choosing $C$ large
enough, depending on $\varepsilon'$. Therefore, we only have to
prove that $P\{V^{\uparrow}(T_{1}^{\uparrow},T_{1}^{\uparrow}+C \log
n) \ge \delta \log n\}= o( n^{-\varepsilon}).$
 Then, applying the strong Markov property at time
$T_{1}^{\uparrow},$ we have to prove that $P\{V^{\uparrow}(0,C \log
n) \ge \delta \log n\}= o( n^{-\varepsilon}).$ Now, by Lemma \ref{lemmeLD} we get
\begin{eqnarray*}
  P\{V^{\uparrow}(0,C \log n) \ge \delta \log n\} &\le&
  (C\log n)^2 \max_{0 \le k \le C \log n} P\{V(k) \ge \delta \log n
\}
\\
   &\le &  (C \log n)^2 \max_{0 \le k \le C \log n} \ee^{-k I\left({\delta
   \log n \over k} \right)}
\\
   &\le & (C \log n)^2 \exp \{- \kappa \delta \log n\}.
\end{eqnarray*}

\noindent Since $\delta > \varepsilon / \kappa,$ this yields
(\ref{Apotfluctueq3}).

To get (\ref{Apotfluctueq2}), observe first that
\begin{eqnarray*}
P\{ -V^{\downarrow}(b,c) \ge \delta \log n\} \le P\{
-V^{\downarrow}(b,T^{\uparrow}_1) \ge \delta \log n \}+P\{
-V^{\downarrow}(T^{\uparrow}_1,c) \ge \delta \log n\}.
\end{eqnarray*}

\noindent The first term on the right-hand side is equal to $P\{
V^{\downarrow}(0,T^{\uparrow}(h_n)) \ge \delta \log n | H_0> h_n
\}.$ Recalling that (\ref{iglehartthm}) implies $P\{ H_0
> h_n\} \le C n^{-(1-\varepsilon)}$ for all large $n$ and
observing the trivial inclusion $\left\{
V^{\downarrow}(0,T^{\uparrow}(h_n)) \ge \delta \log n \, ; \, H_0 >
h_n \right\} \!\subset\! \left\{ T^{\downarrow}(\delta \log
n)\!<T_{h_n}\!< T_{(-\infty,0]} \right\},$ it follows that $P\{
-V^{\downarrow}(b,T^{\uparrow}_1) \ge \delta \log n\}$ is less or
equal than
\begin{eqnarray*} && C
n^{1-\varepsilon} P\{ T^{\downarrow}(\delta \log n)<T_{h_n}<
T_{(-\infty,0]}\} \\& \le & C n^{1-\varepsilon} \sum_{p=\lfloor
\delta \log n \rfloor}^{ \lfloor h_n \rfloor} P\{M_\delta \in
[p,p+1) \, ; \, T^{\downarrow}(\delta \log n)<T_{h_n}<
T_{(-\infty,0]}\},
\end{eqnarray*}

\noindent where $M_\delta:=\max\{V(k); \, 0 \le k \le
T^{\downarrow}(\delta \log n)\}.$ Applying the strong Markov
property at time $T^{\downarrow}(\delta \log n)$ and recalling
(\ref{tailfeller}) we bound the term of the previous sum, for
$\lfloor \delta \log n \rfloor \le p \le \lfloor h_n \rfloor$ and
all large $n,$ by
\begin{eqnarray*} P\{ S \ge p \} \, P\{ S
\ge h_n-(p- \delta \log n)\} \le  C \ee^{-\kappa p} \ee^{-\kappa
(h_n-p+\delta \log n))},
\end{eqnarray*}

\noindent where $S=\sup \{V(k); \, k \ge 0\}.$ Thus, we get $ P\{
-V^{\downarrow}(b,T^{\uparrow}_1) \ge \delta \log n\} \le C \lfloor
h_n \rfloor n^{- \kappa \delta},$ for all large $n,$ which yields $
P\{ -V^{\downarrow}(b,T^{\uparrow}_1) \ge \delta \log n\} =
o(n^{-\varepsilon}),$ $n \to \infty,$ since $\delta> \varepsilon
/\kappa.$ Furthermore, applying the strong Markov property at
$T^{\uparrow}_1,$ we obtain that $P\{
-V^{\downarrow}(T^{\uparrow}_1,c) \ge \delta \log n \}\le P\{
-V^{\downarrow}(0,V_{max}) \ge \delta \log n\}.$ In a similar way we
used before (but easier), we get, by applying the strong Markov
property at $T^{\downarrow}(\delta \log n),$ that $P\{
-V^{\downarrow}(T^{\uparrow}_1,c) \ge \delta \log n\} \le n^{-
\kappa \delta}$ for all large $n.$ Since $\delta> \varepsilon
/\kappa$ this yields (\ref{Apotfluctueq2}).

For (\ref{Apotfluctueq1}), observe first that $((V(k-b)-V(b))_{a \le
k \le b},a,b)$ has the same distribution as $((V(k))_{a^- \le k \le
0},a^-,0)$ under $P\{\cdot |V(k)\ge 0 \, , \, a^- \le k \le 0\},$
where $a^-:=\sup \{ k \le 0 : \, V(k) \ge D_n  \}$. Then, since
$P\{V(k)\ge 0\, , \, k \le 0 \}>0$ and since $(V(-k)\, , \, k \ge
0)$ has the same distribution as $(-V(k)\, , \, k \ge 0),$ we obtain
$$
P\{  V^{\uparrow}(a,b) \ge \delta \log n \} \le C P\{
V^{\uparrow}(0,T_{(-\infty,-D_n]})\ge \delta \log n \}.
$$

\noindent Now, the arguments are the same as in the proof of
(\ref{Apotfluctueq3}). \end{proof}

\medskip

\subsection{Proof of Proposition \ref{P:prop2}}

Recall that we can write $ \tau(d)=\sum_{i=1}^{N} F_i +G,$ where the
$F_i$'s are the successive i.i.d. failures and $G$ the first
success. Then, denoting $F_1$ by $F,$ we have
\begin{eqnarray}
E_{\omega,|a}^{b}[\ee^{- \lambda_n \tau(d)}] &=& E_{\omega,|a}^{b} [
\ee^{- \lambda_n G}] \sum_{k \ge 0} E_{\omega,|a}^{b}[\ee^{-
\lambda_n F}]^k (1-p) p^k
 \nonumber
\\
   &=& E_{\omega,|a}^{b}[\ee^{-
\lambda_n G}] {1-p \over 1-p \, E_{\omega,|a}^{b}[\ee^{- \lambda_n
F}]}. \label{F+S}
\end{eqnarray}

\noindent In order to replace $E_{\omega,|a}^{b}[ \ee^{- \lambda_n
F}]$ by $1-\lambda_n E_{\omega,|a}^{b}[ F],$ we observe that
$1-\lambda_n E_{\omega,|a}^{b}[F] \le E_{\omega,|a}^{b}[\ee^{-
\lambda_n F} ] \le 1-\lambda_n E_{\omega,|a}^{b}[F]+
\frac{\lambda_n^2}{2} E_{\omega,|a}^{b}[ F^2],$ which implies that
$E[{1-p \over 1-p \, E_{\omega,|a}^{b}[ \ee^{- \lambda_n F}]}]$
belongs to
\begin{eqnarray*}
\left[E\bigg[{1-p \over 1-p (1-\lambda_n E_{\omega,|a}^{b}[
F])}\bigg]\, ; \,E\bigg[{1-p \over 1-p (1-\lambda_n
E_{\omega,|a}^{b}[F]+ \frac{\lambda_n^2}{2}
E_{\omega,|a}^{b}[F^2])}\bigg]\right].
\end{eqnarray*}

\noindent Now, we have to bound $\lambda_n E_{\omega,|a}^{b}[F^2]$
from above.  Then, recalling (\ref{expF^2}), which implies $
E_{\omega,|a}^{b}[ F^2]\le 4 ( R^+ + R^-),$ we only have to bound
$R^+$ and $R^-.$ By definition of $R^+,$ we obtain
\begin{eqnarray}
\label{R+majo1}
 R^+ \le(d-b)\left(1+2(d-b)
\ee^{-\widehat{V}^{\downarrow}(b,d)} \right) \left(3(d-b)
 \max_{b \le j \le d} \ee^{-(\widehat{V}(j)-\widehat{V}(b))} \right).
\end{eqnarray}

\noindent Recalling that the estimates (\ref{Aboxlengtheq0})--(\ref{Aboxlengtheq3}) imply that $P\{d-a \ge C'' \log n
\}=o(n^{-\varepsilon})$ and that Lemma \ref{l:Apotfluctu} tells that
$P\{A_5(n)\}=1-o(n^{-\varepsilon}),$ we are interested in the event
$A^{\ddag}(n):=\{d-a \le C'' \log n \} \cap A_5(n),$ whose
probability is greater than $1-o(n^{-\varepsilon})$ for $n$ large
enough. It allows us to sculpt the deep valley $(a,b,c,d)$, such
that we can bound $R^+.$ We are going to show that the fluctuations
of $\widehat{V}$ are, in a sense, related to the fluctuations of $V$
controlled by $A_5(n).$ Indeed, (\ref{Vcheckeq2}) yields
$\widehat{V}^{\downarrow}(b,c)\ge V^{\downarrow}(b,c)\ge - \delta
\log n$ on $A^{\ddag}(n).$ Moreover, (\ref{Vcheckeq0}) together with
(\ref{Vcheckeq1}) imply that $\widehat{V}(y)-\widehat{V}(x)$ is
greater than
\begin{eqnarray*}
 [V(y)-\max_{y \le j \le d-1} V(j)]-[V(x)-\max_{x \le j \le
d-1} V(j)]-O( \log_2 n),
\end{eqnarray*}

\noindent for any $c \le x \le y \le d,$ on $A^{\ddag}(n).$ Since
$V(x)-\max_{x \le j \le d-1} V(j) \le 0$ and $V(y)-\max_{y \le j \le
d-1} V(j)\ge-\delta \log n$ on $A^{\ddag}(n),$ this yields
$\widehat{V}^{\downarrow}(c,d)\ge- \delta \log n-O( \log_2 n).$
Furthermore, since (\ref{Vcheckeq0}) and (\ref{Vcheckeq1}) imply
that $\widehat{V}(c)$ is larger than $\max_{b \le j \le
c}\widehat{V}(j)-O(\log_2 n),$ assembling
$\widehat{V}^{\downarrow}(b,c)\ge - \delta \log n$ with
$\widehat{V}^{\downarrow}(c,d)\ge- \delta \log n-O( \log_2 n)$ yield
\begin{eqnarray}
\label{fluctu(b,d):mino1} \widehat{V}^{\downarrow}(b,d)\ge-
\delta\log n-O( \log_2 n),
\end{eqnarray}

\noindent on $A^{\ddag}(n).$ Therefore, we have, on $A^{\ddag}(n)$
and for all large $n,$
\begin{eqnarray}
\label{R+majo1'} R^+ \le C (\log n)^3 n^{\delta} \max_{b \le j \le
d} \ee^{-(\widehat{V}(j)-\widehat{V}(b))}.
\end{eqnarray}

\noindent Since $\widehat{V}(b)=V(b)$ and (\ref{Vcheckeq1}) implies
$\widehat{V}(x)\ge V(x),$ for all $b \le x \le c$ (in particular
$\widehat{V}(c)\ge V(c)$), it follows from (\ref{fluctu(b,d):mino1})
that
$\widehat{V}(j)-\widehat{V}(b)=(\widehat{V}(j)-\widehat{V}(c))+(\widehat{V}(c)-\widehat{V}(b))
\ge h_n - \delta  \log n -O( \log_2 n),$ which is greater than $0$
for $n$ large enough whenever $ \delta < {(1-\varepsilon) /\kappa}$
(it is possible since $\delta>\varepsilon/\kappa$ and
$0<\varepsilon<1/3$). Therefore, recalling (\ref{R+majo1'}), we
obtain, on $A^{\ddag}(n),$
\begin{eqnarray}
\label{R+majo2}
 R^+ \le C (\log n)^3 n^{\delta}.
\end{eqnarray}

\noindent In a similar way, we prove that $R^- \le C (\log n)^3
n^{\delta},$ on $A^{\ddag}(n),$ which implies that $\lambda_n
E_{\omega,|a}^{b}[ F^2] \le C (\log n)^3
n^{\delta-\frac{1}{\kappa}}.$ Now, observe that, for any $\xi>0,$
$\{\lambda_n E_{\omega,|a}^{b}[ F^2] \le 2(1-\ee^{-\xi})\}$ is
included in $A^{\ddag}(n),$ so that $\lambda_n E_{\omega,|a}^{b}[
F^2] \le 2(1-\ee^{-\xi})E_{\omega,|a}^{b}[ F]$ with probability
larger than $1-o(n^{-\varepsilon}).$ Then, introducing
$$
R'_n(\lambda):=E\bigg[{1 \over 1+\frac{\lambda}{n^{1/\kappa}}
\frac{p}{1-p} E_{\omega,|a}^{b}[F]}\bigg],
$$

\noindent we get, for $n$ large enough,
\begin{eqnarray}
\label{xi0}
 R'_n(\lambda) +o(n^{-\varepsilon}) \le E\bigg[ {1-p \over 1-p \,
E_{\omega,|a}^{b}[ \ee^{- \lambda_n F}]} \bigg]  \le R'_n(\ee^{-\xi}
\lambda) +o(n^{-\varepsilon}).
\end{eqnarray}

In order to bound $E_{\omega,|a}^{b} \left[ \ee^{- \lambda_n G}
\right]$ by below, we observe that $\ee^{-x}\ge 1-x,$ for any $x \ge
0,$ such that $E_{\omega,|a}^{b}[ \ee^{- \lambda_n G}]\ge
1-\lambda_n E_{\omega,|a}^{b}[G].$ Therefore, we only have to bound
$E_{\omega,|a}^{b}[G]$ from above. Recalling (\ref{expS}), we get
$E_{\omega,|a}^{b}[G] \le (d-b)^2 \ee^{\bar{V}^{\uparrow}(b,d)}.$
Now, let us bound $\bar{V}^{\uparrow}(b,d).$ We observe first that
(\ref{Vbareq2}) implies $\bar{V}^{\uparrow}(c,d) \le
V^{\uparrow}(c,d),$ which yields $\bar{V}^{\uparrow}(c,d) \le \delta
\log n$ on $A^{\ddag}(n).$ Moreover, (\ref{Vbareq0}) together with
(\ref{Vbareq1}) imply that $\bar{V}(y)-\bar{V}(x)$ is less or equal
than
\begin{eqnarray*}
 [V(y)-\max_{b \le j \le y} V(j)]-[V(x)-\max_{b \le j \le x} V(j)]+O( \log_2
 n),
\end{eqnarray*}

\noindent for any $b \le x \le y \le c,$ on $A^{\ddag}(n).$ Since
$V(y)-\max_{b \le j \le y} V(j) \le 0$ and $V(x)-\max_{b \le j \le
x} V(j)\ge -\delta \log n$ on $A^{\ddag}(n),$  this yields
$\bar{V}^{\uparrow}(b,c)\le \delta \log n+O( \log_2 n).$
Furthermore, (\ref{Vbareq2}) and the fact that $V(y) \le V(c),$ for
$c\le y \le d,$ imply that $\bar{V}(y) \le \bar{V}(c)$ for $c \le y
\le d.$ Therefore, we have
\begin{eqnarray*}
 \bar{V}^{\uparrow}(b,d)\le \delta\log n+O(
\log_2 n),
\end{eqnarray*}

\noindent on $A^{\ddag}(n).$ It means that $E_{\omega,|a}^{b}[
\ee^{- \lambda_n G}]$ is greater than $1-o(n^{-\varepsilon})$ on
$A^{\ddag}(n)$ whenever $\delta< \frac{1}{\kappa}-\varepsilon,$
which is possible since $\delta>\varepsilon/\kappa$ and
$0<\varepsilon<1/3.$ Therefore, recalling (\ref{xi0}), we obtain
\begin{eqnarray}
\label{xi1}
 R'_n(\lambda) +o(n^{-\varepsilon}) \le
E\big[E_{\omega,|a}^{b} [\ee^{- \lambda_n \tau(d)}] \big]  \le
R'_n(\ee^{-\xi} \lambda) +o(n^{-\varepsilon}).
\end{eqnarray}

\noindent Recalling (\ref{expF}) and (\ref{1-p}), we get
\begin{eqnarray*}
 R_n(\lambda, 2 \widehat{M}_1 (\ee^{H^{(1)}} M_2+\omega_b))\le
R'_n(\lambda)  \le R_n(\lambda, 2 \ee^{H^{(1)}} \widehat{M}_1 M_2),
\end{eqnarray*}

\noindent where $\widehat{M}_1:=\sum_{x=a+1}^{d-1}
\ee^{-(\widehat{V}(x)-\widehat{V}(b))},$ $M_2:=\sum_{x=b}^{d-1}
\ee^{V(x)-V(c)}$ and $ R_n(\lambda,Z)$ is defined in
(\ref{RnlambdaZ}). Furthermore, since $\ee^{H^{(1)}}\ge
n^{\frac{1-\varepsilon}{\kappa}},$ $M_2 \ge 1$ and $\omega_b\le1$ we
obtain that, for any $\xi>0$ and $n$ large enough, $\omega_b \le
(\ee^{\xi}-1)\ee^{H^{(1)}}M_2.$ Therefore, we have for all large
$n,$
\begin{eqnarray}
\label{xi2}
 R_n(\ee^{\xi}\lambda, 2  \ee^{H^{(1)}}\widehat{M}_1 M_2)\le
R'_n(\lambda)  \le R_n(\lambda, 2 \ee^{H^{(1)}} \widehat{M}_1 M_2).
\end{eqnarray}

\noindent Now, assembling (\ref{xi1}) and (\ref{xi2}) concludes the
proof of Proposition \ref{P:prop2}. \hfill$\Box$

\section{Back to canonical meanders}
\label{sec:part3}

Recall that $S=\max\{ V(k)\, ; \, k\ge 0\}$ and let us set $H:=\max\{V(k) \, ; \,
0\le k\le T_{\r_-}\}=H_0,$ $T_S:=\inf\{k\ge 0: V(k)=S\}.$
Moreover, we define $\mathcal{I}_n:=\{H=S \ge h_n\}\cap \{V(k) \ge
0\, , \, \forall \, k\le 0\},$ and introduce the random variable
$Z:=e^S M_1^+ M_2^+,$ where $M_1^+:= \sum_{k=a^-}^{T_{h_n/2}}
e^{-V(k)}$ and $M_2^+:=\sum_{k=0}^{d^+} e^{V(k)-S},$ with $a^-=\sup
\{ k \le 0 : \, V(k) \ge D_n  \}$ and $d^+:=\inf\{ k \ge e_1 : \,
V(k)-V(e_1) \le -D_n  \}.$ Then, denoting
$$
\mathcal{R}_n(\lambda):=E\bigg[{1 \over 1+n^{-\frac{1}{\kappa}}2
\lambda Z}|\mathcal{I}_n \bigg],
$$

\noindent we get the following result.
\medskip

\begin{proposition}
 \label{P:prop3} For any $\xi>0,$ we have, for $n$ large enough,
 \begin{eqnarray*}
\mathcal{R}_n(\ee^{\xi} \lambda) +o(n^{-\varepsilon}) \le
R_n(\lambda, 2 \ee^{H^{(1)}} \widehat{M}_1 M_2)\le
\mathcal{R}_n(\ee^{-\xi}\lambda) +o(n^{-\varepsilon}).
\end{eqnarray*}
\end{proposition}

\medskip

\begin{proof} {\it Step 1: we replace $\widehat{M}_1$ by
$\widehat{M}_1^T.$}

 Recall that $A^{\ddag}(n)=\{d-a \le C'' \log n \} \cap
A_5(n)$ and that $P\{A^{\ddag}(n)\} \ge 1-o(n^{-\varepsilon}),$ for
all large $n.$ Now, let us introduce $T(\frac{h_n}{2}):=\inf \{k \ge
b :V(k)-V(b)\ge h_n/2 \}$ and $\widehat{M}_1^T:=
\sum_{k=a+1}^{T(\frac{h_n}{2})}e^{-(\widehat{V}(k)-\widehat{V}(b))}.$
Recalling (\ref{fluctu(b,d):mino1}), we observe that $\widehat{M}_1
\le \widehat{M}_1^T + C'' \log n \ee^{-\frac{h_n}{2}+\delta \log n}$
on $A^{\ddag}(n).$ This implies that, for any $\xi>0,$ we have
$\widehat{M}_1 - \widehat{M}_1^T\le(\ee^{\xi}-1)\widehat{M}_1^T$ for
all large $n,$ whenever $\delta<\frac{1-\varepsilon}{2 \kappa},$
which is possible since $\delta>\varepsilon/\kappa$ and
$0<\varepsilon<1/3.$ Therefore, we obtain, for $n$ large enough,
\begin{eqnarray*}
R_n(\ee^{\xi} \lambda, 2 \ee^{H^{(1)}} \widehat{M}^T_1
M_2)+o(n^{-\varepsilon}) \le R_n(\lambda, 2 \ee^{H^{(1)}}
\widehat{M}_1 M_2) \le R_n(\lambda, 2 \ee^{H^{(1)}} \widehat{M}^T_1
M_2).
\end{eqnarray*}

\medskip

{\it Step 2: we replace $\widehat{M}_1^T$ by $M_1^T.$}

Let us denote $M_1^T:=\sum_{k=a+1}^{T(\frac{h_n}{2})}
e^{-(V(k)-V(b))}.$ Since $T(\frac{h_n}{2}) \le c,$ (\ref{Vcheckeq2})
implies that $\widehat{M}_1^T \le M_1^T.$  Observe that
(\ref{Vcheckeq0}) with (\ref{Vcheckeq1}) imply that
$\widehat{V}(y)-\widehat{V}(b)- (V(y)-V(b))$ is less or equal than
\begin{eqnarray*}
&&\log \bigg(\frac{\sum_{j=b}^{d-1} \ee^{V(j)}}{\sum_{j=y}^{d-1}
\ee^{V(j)}} \frac{\sum_{j=b+1}^{d-1} \ee^{V(j)}}{\sum_{j=y+1}^{d-1}
\ee^{V(j)}} \bigg) \le \frac{\sum_{j=b}^{y-1}
\ee^{V(j)}}{\sum_{j=y}^{d-1} \ee^{V(j)}} + \frac{\sum_{j=b+1}^{y}
\ee^{V(j)}}{\sum_{j=y+1}^{d-1} \ee^{V(j)}},
\end{eqnarray*}

\noindent for any $b \le y \le d.$ Therefore, on $A^{\ddag}(n),$ we
obtain $\widehat{V}(y)-\widehat{V}(b) \le (V(y)-V(b))+ C \log n
\ee^{-\frac{ h_n}{2}}$ for any $b \le y \le T(\frac{h_n}{2}),$ which
yields  $\widehat{M}_1^T \ge \exp\{C \log n \, \ee^{-\frac{
h_n}{2}}\} M_1^T.$ Then, for any $\xi>0,$ we obtain that
$\widehat{M}_1^T \ge \ee^{-\xi} M_1^T,$ on $A^{\ddag}(n)$ and for
all large $n.$ This implies
\begin{eqnarray*}
R_n(\lambda, 2 \ee^{H^{(1)}} M^T_1 M_2) \le R_n(\lambda, 2
\ee^{H^{(1)}} \widehat{M}^T_1 M_2) \le R_n(\ee^{-\xi} \lambda, 2
\ee^{H^{(1)}}M^T_1 M_2)+o(n^{-\varepsilon}).
\end{eqnarray*}

Now, assembling Step $1$ and Step $2,$ we get that, for any $\xi>0$
and $n$ large enough, $R_n(\lambda, 2 \ee^{H^{(1)}} \widehat{M}_1
M_2)$ belongs to
\begin{eqnarray}
\label{xi3} \left[R_n(\ee^{\xi} \lambda, 2 \ee^{H^{(1)}} M^T_1
M_2)+o(n^{-\varepsilon}) \, ; \,  R_n(\ee^{-\xi} \lambda, 2
\ee^{H^{(1)}}M^T_1 M_2)+o(n^{-\varepsilon})\right].
\end{eqnarray}

 {\it Step 3: the ``good '' conditioning.}

Let us first observe that $((V(k-b)-V(b))_{a\le k \le d},a,b,c,d)$
has the same law as $((V(k))_{a^-\le k \le d^+}, a^-,0, T_H, d^+)$
under $P\{\cdot | \mathcal{I}'_n\},$ where $\mathcal{I}'_n:=\{H \ge
h_n\, ; \, V^\uparrow(a^-,0) \le h_n \, ; \, V(k) \ge 0 \, , \, a^-
\le k \le 0\}.$ Moreover, we easily obtain that $P\{\{ V(k) \ge 0 \,
, \, a^- \le  k \le 0\} \setminus \{V(k) \ge 0 \, , \, k \le 0 \}
\}=O(n^{-(1+\kappa)})=o(n^{-\varepsilon}),$ that $P\{\{ H \ge h_n\}
\setminus \{H=S\} \}=O(n^{-2(1-\varepsilon)})=o(n^{-\varepsilon})$
and that $P\{V^\downarrow(a^-,0) > h_n\}\le P\{V^\downarrow(a^-,0) >
\delta \log n\}=o(n^{-\varepsilon}),$ with the same arguments as in
the proof of Lemma \ref{l:Apotfluctu}. Therefore, we have $P\{
\mathcal{I}'_n \bigtriangleup \mathcal{I}_n \}=o(n^{-\varepsilon}).$
Since $0 \le R_n(\lambda,Y) \le 1,$ for any $\lambda>0$ and any
positive random variable $Y,$ this yields
\begin{eqnarray}
\label{goodlaw} R_n(\lambda, 2 \ee^{H^{(1)}} M^T_1
M_2)=\mathcal{R}_n(\lambda)+o(n^{-\varepsilon}).
\end{eqnarray}

\noindent Combining (\ref{xi3}) and (\ref{goodlaw}) together concludes
the proof of Proposition \ref{P:prop3}.\end{proof}

\section{Proof of Theorem \ref{t:main}}
\label{sec:part4}

Observe first that $\mathcal{R}_n(\lambda)$ can be written
$$
\mathcal{R}_n(\lambda)=1-E\bigg[1-{1 \over 1+2 \lambda_n
Z}|\mathcal{I}_n \bigg].
$$

\noindent Then, we can use Corollary A.1 and Remark A.1 in
\cite{enriquez-sabot-zindy+}, that together imply
$$
E\bigg[1-{1 \over 1+ 2 \lambda_n Z}
 \, \big| \,   \mathcal{I}_n \bigg] \sim 2^\kappa \frac{\pi \kappa}{\sin(\pi
 \kappa)} \ \frac{E[M^\kappa]^2 C_I}{n P\{H\ge h_n\}} \  \lambda^\kappa, \qquad n \to
 \infty.
$$
where the random variable $M$ defined by 
\begin{equation}
\label{eq:defM}
M:=\sum_{k<0} \ee^{-V'_k}+ \sum_{k\ge 0}
\ee^{-V''_k},
\end{equation}
where $(V'_k)_{k < 0}$ is distributed as the potential under $P\{\cdot|V_k \ge 0, \,
\forall k <0\}$ while  $(V''_k)_{k \ge 0}$ is independent of $(V'_k)_{k < 0}$ and is distributed as the potential under
$\tilde{P}\{\cdot|V_k > 0, \, \forall k > 0\}.$

\noindent Therefore, combining together the results of Proposition \ref{P:prop1},
Proposition \ref{P:prop2}, Proposition \ref{P:prop3} and recalling
that $q_n:=P\{H\ge h_n\},$ we get that, for any $\xi>0,$
\begin{eqnarray*}
\liminf_{n \to \infty} \e[\ee^{- \lambda_n \, \tau(e_n)}] &\ge&
\exp\Big\{-\Big(2^\kappa \frac{\pi \kappa}{\sin(\pi
 \kappa)}E[M^\kappa]^2
C_I\Big)(\ee^\xi\lambda)^\kappa\Big\},
\\
\limsup_{n \to \infty} \e[\ee^{- \lambda_n \, \tau(e_n)}] &\le&
\exp\Big\{-\Big(2^\kappa \frac{\pi \kappa}{\sin(\pi
 \kappa)}E[M^\kappa]^2
C_I\Big)(\ee^{-\xi}\lambda)^\kappa \Big\}.
\end{eqnarray*}

\noindent Since this result holds for any $\xi>0,$ we get,
\begin{eqnarray*}
\lim_{n \to \infty} \e[\ee^{- \lambda_n \,
\tau(e_n)}]=\exp\Big\{-\Big(2^\kappa \frac{\pi \kappa}{\sin(\pi
\kappa)}E[M^\kappa]^2 C_I\Big)\lambda^\kappa\Big\}.
\end{eqnarray*}

\noindent 

Now, one can
be tempted to express the functional $E[M^\kappa]$ in terms of the
more usual constant $C_K,$ see (\ref{tailkesten}). This is the
content of Theorem $2.1$ in \cite{enriquez-sabot-zindy+}, which
yields
$$ C_K=E[M^\kappa]C_F  =E[M^\kappa] {(1-E[\ee^{\kappa V(e_1)}])\over\kappa E[\rho_0^\kappa\log\rho_0]E[e_1]}.$$

\noindent Therefore, the Laplace transform of $ n^{-1/\kappa}
\tau(e_n)$ is
$$\begin{array}{rl}\ds\E[\ee^{-{\lambda \over n^{1/\kappa}} \tau(e_n)}]&=\ds \exp\Big\{-\Big(2^\kappa
 {\pi\kappa\over\sin(\pi\kappa)}{C_K^2 C_I\over C_F^2}\Big)\lambda^\kappa\Big\}+o(1)\\&
=\ds \exp\Big\{-\Big(2^\kappa
{\pi\kappa^2\over\sin(\pi\kappa)}C_K^2E[\rho_0^\kappa\log\rho_0]E[e_1]\Big)\lambda^\kappa\Big\}+o(1).\end{array}$$
Finally, since, by the law of large numbers, $e_n/n$ converges almost surely
to $E[e_1]$, we conclude that
$$\E[\ee^{-{\lambda \over n^{1/\kappa}} \tau(n)}]
=\ds \exp\Big\{-\Big(2^\kappa
{\pi\kappa^2\over\sin(\pi\kappa)}C_K^2E[\rho_0^\kappa\log\rho_0]\Big)\lambda^\kappa\Big\}+o(1).$$
Hence, we obtain that the limit is the positive stable law with
index $\kappa$ and parameter $2^\kappa
{\pi\kappa^2\over\sin(\pi\kappa)}C_K^2E[\rho_0^\kappa\log\rho_0].$
\hfill$\Box$
\medskip

We can easily see that we can deduce from this proof the
asymptotic of the Laplace transform of the time needed to cross
the first $*$-valley.
\begin{corollary}
We have
$$
\E\left[ 1-e^{-{\lambda\over n^{1/\kappa}}\tau_1^*}\right] \sim
2^\kappa{\pi\kappa\over \sin(\pi\kappa)}{C_U\over n P(H\ge
h_n)}\lambda^\kappa,
$$
where $C_U=C_I E[M^\kappa]$ is the constant which appears in the
tail estimate of $Z$, in \cite{enriquez-sabot-zindy+}.
\end{corollary}
\begin{remark}
This result would hold for a different choice of $h_n$. Indeed,
from the proof of Proposition \ref{P:prop2} and Proposition \ref{P:prop3} and
from Corollary A.1 of \cite{enriquez-sabot-zindy+}, we see that
the result holds for any choice of $h_n$ such that
$e^{h_n}=o(n^{1\over \kappa})$ and $h_n\ge n^{1-\epsilon\over
\kappa}$ for some $0<\epsilon<1/3$ (this last condition comes
from the technical assumption in (4.1) which is needed in the proof
Proposition \ref{P:prop2}, see (\ref{R+majo2})).
\end{remark}

\section{Proof of Corollary 1}

We are in the case  when
the law of the environment satisfies $$\omega_{1}(\d
x)=\frac{1}{B(\alpha,\beta)}x^{\alpha-1}(1-x)^{\beta-1}{\bf
1}_{[0,1]}(x) \d x,$$ with $\alpha, \beta >0$ and
$B(\alpha,\beta):=\int_{0}^1 x^{\alpha-1}(1-x)^{\beta-1} \d x,$.
 The assumption of Theorem
\ref{t:main} corresponds to the case where $0< \alpha-\beta <1$ and
an easy computation leads to $\kappa=\alpha-\beta.$
 Now, a classical argument of
derivation under the sign integral shows that
$$ E[\rho_0^\kappa\log\rho_0]=\psi(\a)-\psi(\b),$$
where $\psi$ denotes the classical Digamma function
$\psi(z):=(\log\Gamma)'(z)={\Gamma'(z)\over\Gamma(z)}$.
Furthermore, a work of Chamayou and Letac \cite{chamayou-letac}
shows that $C_K$ can be made explicit. Indeed, with the notations of
\cite{chamayou-letac}, $\rho_0$ follows the law
$\beta^{(2)}_{p,q}(dx):={1\over
B(p,q)}x^{p-1}(1+x)^{-p-q}1_{\R_+}(x)dx$ with $p=\b$ and $q=\a$.
Then, Example 9 of \cite{chamayou-letac} says that $\sum_{k\geq1}
\ee^{V(k)}$ follows the  law of $\beta^{(2)}_{\b,\a-\b}$ having
density ${1\over B(\a,\b)}x^{\b-1}(1+x)^{-\a}1_{\R_+}(x)$. But we
have $\beta^{(2)}_{\b,\a-\b}([t,+\infty[)\sim {1\over
(\a-\b)B(\a,\b)}{1\over t^{\a-\b}},$ $t \to \infty.$ Hence,
$C_K={1\over (\a-\b)B(\a,\b)}$.

\section{Toward the case $\kappa=1$}
\label{sec:extension} We intend to treat the critical case
$\kappa=1$ between the transient ballistic and sub-ballistic cases.
This case turns out to be more delicate. Indeed, Lemma
\ref{lemmaenv} is replaced by a weaker statement, which says that
$\tau(e_n)$ reduces to the time spent by the walker to climb
excursions which are higher than $\alpha \log n$ for $\alpha$
arbitrarily small. Due to this reduced height, the new ``high''
excursions are much more numerous and are not anymore well
separated. The definition of the valleys should then be adapted as
well as the ``linearization'' argument, which is more difficult to
carry out. Moreover, a result of Goldie \cite{goldie} gives an
explicit formula for the Kesten's renewal constant, namely
$C_K=\frac{1}{E[\rho_0 \log \rho_0]}.$ As a result, we should
obtain,  as a consequence  of a fluctuation result, the following result, which takes a remarkably
simple form: $X_n/({n\over \log n})$ converges in probability to
$E[\rho_0 \log \rho_0]/2.$

\bigskip
\bigskip
\noindent {\bf Acknowledgements}
 Many thanks are due to an anonymous referee for
careful reading of the original manuscript and for helpful comments.

\bigskip

\end{document}